\documentclass[aos,seceqn,nameyear,dvips]{arximspdf}
\usepackage{graphics}
\doi{10.1214/10-AOS808}
\volume{38}
\issue{5}
\pubyear{2010}
\firstpage{2857}
\lastpage{2883}

\makeatletter

\newtheorem{proposition}{Proposition}[section]
\newtheorem{theorem}[proposition]{Theorem}

\newcommand{\BF}{\mathrm{BF}}
\newcommand{\SBF}{\mathrm{SBF}}

  \let\sv@tabnotetext\tabnotetext
  \let\sv@tabnotemark@fmt\tabnotemark@fmt
   \long\def\legend#1{{\let\tabnote@indent\leavevmode\sv@tabnotetext[]{}{#1}}}

\makeatother

\begin{document}
\begin{frontmatter}

\title{Backfitting and smooth backfitting for additive quantile models}
\runtitle{Additive quantile models}

\begin{aug}
\author[A]{\fnms{Young Kyung} \snm{Lee}\thanksref{t1}\ead[label=e1]{youngklee@kangwon.ac.kr}},
\author[B]{\fnms{Enno} \snm{Mammen}\thanksref{t2}\ead[label=e2]{emammen@rumms.uni-mannheim.de}} and
\author[C]{\fnms{Byeong U.} \snm{Park}\corref{}\thanksref{t3}\ead[label=e3]{bupark@stats.snu.ac.kr}}
\runauthor{Y. K. Lee, E. Mammen and B. U. Park}
\affiliation{Kangwon National University, University of Mannheim and
Seoul~National~University}
\address[A]{Y. K. Lee\\
Department of Statistics\\
Kangwon National University\\
Chuncheon 200-701\\
Korea\\
\printead{e1}} %adresu isvedimo komanda gale!
\address[B]{E. Mammen\\
Department of Economics \\
University of Mannheim \\
68131 Mannheim, L7, 3-5\\
Germany\\
\printead{e2}}
\address[C]{B. U. Park\\
Department of Statistics\\
Seoul National University\\
Seoul 151-747\\
Korea\\
\printead{e3}}
\end{aug}

\thankstext{t1}{Supported by Basic Science Research
Program through the National Research Foundation of Korea (NRF) funded
by the Ministry of Education, Science and Technology (2009-0058380).}

\thankstext{t2}{Supported by DFG-project MA 1026/9-2 of the
Deutsche Forschungsgemeinschaft.}

\thankstext{t3}{Supported by Mid-career Researcher Program through NRF grant
funded by the MEST (No. 2010-0017437).}

% HISTORY:
\received{\smonth{10} \syear{2009}}
\revised{\smonth{1} \syear{2010}}

% ABSTRACT
%
\begin{abstract}
In this paper, we study the ordinary backfitting and smooth backfitting
as methods of fitting additive quantile models. We show that these
backfitting quantile estimators are asymptotically equivalent to the
corresponding backfitting estimators of the additive components in a
specially-designed additive mean regression model. This implies that
the theoretical properties of the backfitting quantile estimators are
not unlike those of backfitting mean regression estimators. We also
assess the finite sample properties of the two backfitting quantile
estimators.
\end{abstract}

% KEYWORDS
%

\end{frontmatter}

%%%%%%%%% Section 1
%s1 ###
\section{Introduction}
\label{sec:1}

Nonparametric additive models are powerful techniques for
high-dimensional data. They enable us to avoid the curse of
dimensionality and estimate the unknown functions in high-dimensional
settings at the same accuracy as in univariate cases. In the mean
regression setting, there have been many proposals for fitting additive
models. These include the ordinary backfitting procedure of
\citet{BHT89}, whose theoretical properties were studied
later by \citet{OR97} and \citet{O00}, the marginal
integration technique of \citet{LN95}, and the smooth
backfitting of \citet{MLN99}, \citet{MP06} and
\citet{YPM08}. It is widely accepted that the
marginal integration method still suffers from the curse of
dimensionality since it does not produce rate-optimal estimates unless
smoothness of the regression function increases with the number of
additive components. On the contrary, the ordinary backfitting and
smooth backfitting are known to achieve the univariate optimal rate of
convergence under certain regularity conditions.

In this paper, we are concerned with nonparametric estimation of
additive conditional quantile functions. Conditional quantile
estimation is also a very useful tool for exploring the structure of
the conditional distribution of a response given a predictor. A
collection of conditional quantiles, when graphed, give a picture of
the entire conditional distribution. It can be used directly to
construct conditional prediction intervals. Also, it may be a basis for
verifying the presence of conditional heteroscedasticity;
see \citet{F04}, for example. Various other applications of
conditional quantile
estimation may be found in \citet{YLS03}. In the
nonadditive setting, there have been many proposals for this problem,
which include the work by \citet{JH90}, \citet{C91},
\citet{YJ98}
and \citet{LLP06}. There have been also
some proposals for additive quantile regression. \citet{FG96}
provided a direct extension of the ordinary backfitting method to
quantile regression, but without discussing its statistical properties.
\citet{LY04} gave a heuristic discussion of the asymptotic limit of
a backfitting local linear quantile estimator. \citet{HL05}
studied an extension of the two-stage procedure of
\citet{HM04} to quantile regression. Their estimator is a one-step kernel
smoothing iteration of an orthogonal series estimator.

The main theme of this paper is to discuss the statistical properties
of the ordinary and smooth backfitting methods in additive quantile
regression. The methods are difficult to analyze since there exists no
explicit definition for the ordinary backfitting estimator and, for
both estimators, the objective functions defining the estimators are
not differentiable. We borrow empirical process techniques to tackle
the problem. In particular, we devise a theoretical mean regression
model by using a Bahadur representation for the sample
quantiles. We show that the least squares ordinary and smooth
backfitting estimators in this theoretical mean regression model are
asymptotically equivalent to the corresponding quantile estimators in
the original model. This makes the theoretical properties of the two
backfitting quantile estimators well understood from the existing
theory for the corresponding least squares backfitting mean regression
estimators. The theory was confirmed by a simulation study. Also, it
was observed in the simulation study that the smooth backfitting
estimator outperformed the ordinary backfitting estimator in additive
quantile regression.

The paper is organized as follows. In the next section, the ordinary
and smooth backfitting methods for additive quantile regression are
introduced and their theoretical properties are provided. In Section
\ref{sec:num}, some computational aspects of the smooth backfitting method are
discussed. The simulation results for the finite sample properties of
the two backfitting methods are presented in Section \ref{sec:sim}. Technical
details are given in Section \ref{sec5}.

%%%%%%%%% Section 2
%s2 ###
\section{Main results}
\label{sec:2}

It is assumed for one-dimensional response variables $Y^1$,\break $\ldots,
Y^n$ that
%
%e2.1 ###
%
\begin{equation}\label{eq:2.1}
Y^i = m_0 +m_1(X_1^i) + \cdots+ m_d(X_d^i) +\varepsilon^i,\qquad
1 \le i \le n.
\end{equation}
Here, $\varepsilon^i$ are error variables, $m_1, \ldots,m_d$ are
unknown functions from $\mathbb{R}$ to $\mathbb{R}$ satisfying $\int
m_j(x_j) w_j(x_j) \,dx_j=0$ for some weight functions $w_j$, $m_0$ is an
unknown constant, and $X^i=(X_1^i,\ldots,X_d^i)$ are random design
points in $\mathbb{R}^d$. Throughout the paper, we assume that $(X^i,
\varepsilon^i)$ are i.i.d. and that $X_j^i$ takes its values in a
bounded interval $I_{j}$. Furthermore, it is assumed that the
conditional $\alpha$-quantile of $\varepsilon^i$ given $X^i$ equals
zero. This model excludes interesting auto-regression models, but it
simplifies our asymptotic analysis. We expect that our results can be extended
to dependent observations under mixing conditions.

The ordinary backfitting estimator is based on an iterative algorithm.
The estimate of $m_j$ is updated by the following equation:
%
%e2.2 ###
%
\begin{eqnarray}\label{it:BF}
\hat m^{\BF}_j(x_j) &=& \mathop{\arg\min}_{\theta\in\Theta}
\sum_{i=1}^{n}\tau_{\alpha} \Biggl(Y^{i}-\theta-\hat m^{\BF}_0-
\sum_{\ell=1, \neq j}^d \hat m^{\BF}_\ell(X_\ell^i)
\Biggr)\nonumber\\[-8pt]\\[-8pt]
&&\hspace*{46.7pt}{} \times K_{j,h_j}(x_j,X_j^{i}). \nonumber
\end{eqnarray}
Here, $\tau_\alpha$ is the so called ``check function'' defined by
$\tau_{\alpha}(u)=u\{\alpha-I(u<0)\}$, and $K_{j,g}$ are kernel
functions with bandwidth $g$; see the assumptions below. To simplify
the mathematical argumentation, the minimization in (\ref{it:BF}) runs
over a compact set $\Theta$. It is assumed that all values of the
function $m_j$ lie in the interior of $\Theta$. As in the case of mean
regression, the ordinary backfitting estimator is not defined as a
solution of a global minimization problem.

The smooth backfitting estimator is also based on an iterative
algorithm. The estimate of $m_j$ is updated by the following integral equation:
%
%e2.3 ###
%
\begin{eqnarray}\label{it:SBF}\qquad
\hat m^{\SBF}_j(x_j) &=& \mathop{\arg\min}_{\theta\in\Theta}
\sum_{i=1}^{n}\int\tau_{\alpha} \Biggl(Y^{i}-\theta-\hat
m^{\SBF}_0-\sum_{\ell=1, \neq j}^d \hat m^{\SBF}_\ell(x_\ell)
\Biggr)\nonumber\\[-8pt]\\[-8pt]
&&\hspace*{58.3pt}{} \times \prod_{\ell=1, \neq j}
K_{\ell,h_\ell}(x_\ell,X_\ell^{i}) \,dx_\ell\cdot
K_{j,h_j}(x_j,X_j^{i}), \nonumber
\end{eqnarray}
where the integration is over the support of $(X_1^i, \ldots,
X_{j-1}^i, X_{j+1}^i, \ldots, X_d^i)$. This is an iterative scheme for
obtaining $\hat{m}^{\SBF}_{j}, j=0, 1, \ldots, d$, which minimize
%
%e2.4 ###
%
\begin{eqnarray} \label{minprob}
&&\sum_{i=1}^{n} \int\tau_\alpha\Biggl( Y^{i} - \hat{m}^{\SBF}_{0} -
\sum_{j=1}^{d} \hat{m}^{\SBF}_{j} (x_{j}) \Biggr)\nonumber\\[-8pt]\\[-8pt]
&&\qquad\hspace*{1.8pt}{}\times K_{1,h_1}(x_1,X_1^{i})
\cdots K_{d,h_d}(x_d,X_d^{i}) \,dx_1 \cdots dx_d,\nonumber
\end{eqnarray}
where the integration is over the support of $X^i$. The minimizations
or iterations are done under the constraints
%
%e2.5 ###
%
\begin{equation} \label{constr0}
\int_{I_j} \hat{m}^{l}_j (x_j) w_j (x_j) \,dx_j
=0,\qquad j=1,\ldots, d \mbox{ and } l= \BF, \SBF
\end{equation}
for some weight functions $w_j$. One may take unknown weight functions
such as the marginal densities of $X_j$ and use consistent estimators
of them as the weight functions $w_j$ in the integrals (\ref
{constr0}). But this would lead to more complicated bias calculation.

We compare our model (\ref{eq:2.1}) with the following theoretical
model. For $i=1, \ldots,n$, let $Z^1, \ldots, Z^n$ be one-dimensional
variables such that
%
%e2.6 ###
%
\begin{equation}\label{eq:2.2}
Z^i = m_0 +m_1(X_1^i) + \cdots+ m_d(X_d^i) +\eta^i.
\end{equation}
Here, the constant $m_0$, the functions $m_1, \ldots,m_d$ and the
covariates $X_1^i,\ldots,X_d^i$ are those in (\ref{eq:2.1}). The error
variables $\eta^i$ are defined by
\[
\eta^i = - \frac{I(\varepsilon^i \leq0) - \alpha}{f_{\varepsilon
|X}(0|X^i)},
\]
where $f_{\varepsilon|X}$ is the conditional density of $\varepsilon$
given $X$. This definition is motivated from the Bahadur representation
of sample quantiles [\citet{B66}]. For an independent sample of
$\varepsilon^1,\ldots,\varepsilon^n$ with densities $f_i$ and
$\alpha$-quantiles being equal to~0, the Bahadur expansion states that
the $\alpha$th sample quantile $\hat\theta_{\alpha}$ of
$\varepsilon^1,\ldots, \varepsilon^n$ is asymptotically equivalent
to the weighted average
\[
{\sum_{i=1}^n f_i(0) \eta^i \over\sum_{i=1}^n f_i(0)},
\]
where $\eta^i = - \{I(\varepsilon^i \leq0) - \alpha\} f_i(0)^{-1}$.
Thus, the estimator $\hat\theta_{\alpha}$ is asymptotically
equivalent to the minimizer of
\[
\theta\to\sum_{i=1}^n f_i(0) (\eta^i - \theta)^2.
\]

This consideration suggests that the ordinary and smooth backfitting
estimators defined at (\ref{it:BF}) and (\ref{it:SBF}), respectively,
may be approximated well by the corresponding weighted local least
squares estimators in the model (\ref{eq:2.2}). Note that the model
(\ref{eq:2.2}) is an additive model with errors
$\eta^i$ having conditional mean zero given the covariates $X^i$.
Thus, the weighted ordinary backfitting estimators $\hat m^{*,\BF}_j$ in
this model are defined by the following iterations:
%
%e2.7 ###
%
\begin{eqnarray}\label{it:TBF}\hspace*{19pt}
\hat m^{*,\BF}_j(x_j) &=& \mathop{\arg\min}_{\theta\in\Theta}
\sum_{i=1}^{n} \Biggl\{Z^{i}-\theta- \hat m^{*,\BF}_0 -
\sum_{\ell=1, \neq j}^d \hat m^{*,\BF}_\ell(X_\ell^i) \Biggr\}^2
\nonumber\\
&&\hspace*{46.4pt}{} \times f_{\varepsilon|X}(0|X^i)
K_{j,h_j}(x_j,X_j^{i}) \nonumber\\[-8pt]\\[-8pt]
&=& \sum_{i=1}^{n} \Biggl\{Z^{i}-\hat m^{*,\BF}_0 - \sum_{\ell=1,
\neq j}^d \hat m^{*,\BF}_\ell(X_\ell^i) \Biggr\}
f_{\varepsilon|X}(0|X^i) K_{j,h_j}(x_j,X_j^{i})\nonumber\\
&&\hspace*{12.5pt}{} \times\Biggl\{\sum_{i=1}^{n}
f_{\varepsilon|X}(0|X^i)K_{j,h_j}(x_j,X_j^{i}) \Biggr\}^{-1}.\nonumber
\end{eqnarray}
Also, the weighted smooth backfitting estimators $\hat m^{*,\SBF}_j$ in
the model (\ref{eq:2.2}) are defined by
%
%e2.8 ###
%
\begin{eqnarray}\label{it:TSBF}
\hat m^{*,\SBF}_j(x_j) &=& \tilde m^{*,\SBF}_j(x_j)
-\hat m^{*,\SBF}_0 \nonumber\\[-8pt]\\[-8pt]
&&{} - \sum_{\ell=1, \ne j}^d \int{\hat
m^{*,\SBF}_\ell(x_\ell) {\hat f^w_{X_j,X_\ell}(x_j,x_\ell) \over
\hat f^w_{X_j}(x_j)} \,dx_\ell},\nonumber
\end{eqnarray}
where
\begin{eqnarray*}
\tilde m^{*,\SBF}_j(x_j)&=& n^{-1}\sum_{i=1}^{n}Z^{i}
f_{\varepsilon|X}(0|X^i) K_{j,h_j}(x_j,X_j^{i}) \hat
f^w_{X_j}(x_j)^{-1},\\
\hat f^w_{X_j}(x_j)&=& n^{-1}\sum_{i=1}^{n}f_{\varepsilon|X}(0|X^i)
K_{j,h_j}(x_j,X_j^{i}),\\
\hat f^w_{X_j,X_\ell}(x_j,x_\ell)&=& n^{-1}\sum_{i=1}^{n}
f_{\varepsilon|X}(0|X^i)
K_{j,h_j}(x_j,X_j^{i})K_{\ell,h_\ell}(x_\ell,X_\ell^{i})
\end{eqnarray*}
are weighted modifications of the marginal Nadaraya--Watson estimator
and the kernel estimators of the one- and two-dimensional marginal
densities of $X$, respectively. The latter two are in fact kernel
estimators of
\begin{eqnarray*}
f_{X_j}^w (x_j) &=& \int f_{\varepsilon|X}(0|x)f_X(x) \,dx_{-j} =
f_{\varepsilon,X_j}(0,x_j),\\
f_{X_j,X_\ell}^w (x_j,x_\ell) &=& \int f_{\varepsilon|X}(0|x)f_X(x)
\,dx_{-(j,\ell)}= f_{\varepsilon,X_j,X_\ell}(0,x_j,x_\ell),
\end{eqnarray*}
respectively, where $x_{-j}=(x_1, \ldots, x_{j-1}, x_{j+1}, \ldots,
x_d)^\top$ and $x_{-(j,\ell)}$ is a vector that has elements $x_l$
with $1 \leq l \leq d$ and $l\not=j,\ell$.

Our first result (Proposition \ref{TH:2.1}) shows that each application of the
updating equations (\ref{it:TBF}) and (\ref{it:TSBF}) in the
theoretical model (\ref{eq:2.2}), respectively, lead to asymptotically
equivalent results with those at (\ref{it:BF}) and (\ref{it:SBF}) in
the original model~(\ref{eq:2.1}). In the next step, we will apply
Proposition \ref{TH:2.1} for iterative applications of the backfitting updates.
We will show
that the asymptotic equivalence remains to hold for iterative
applications of the backfitting procedures as long as the number of
iterations is small enough. By extending the results for backfitting
and smooth backfitting estimators in mean regression, we will use this
fact to get our main result (Theorem \ref{TH:2.2}). The latter states an
asymptotic normality result for the ordinary and smooth backfitting
quantile estimators in additive models. Its proof is based on an
argument that carries an asymptotic normality result in mean regression
over to quantile regression.

We now introduce assumptions that guarantee asymptotic equivalence
between the mean and the quantile backfitting estimators after one
cycle of update. Further assumptions that are needed for iterative
updates will be given after Proposition~\ref{TH:2.1}. For simplicity, we state
Proposition \ref{TH:2.1} and its conditions only for the updates of the first
additive component. In abuse of notation, we denote the estimators of
the components $m_j, 2 \le j \le d$, at the preceding iteration
step, by $\hat m_2^{l}, \ldots, \hat m_d^{l}$, where $l$ stands for
$\BF$, $ \SBF$, $*,\BF$ or $*,\SBF$. The updates of the first component
that are obtained by plugging these estimators into (\ref{it:BF}),
(\ref{it:SBF}), (\ref{it:TBF}) and (\ref{it:TSBF}), respectively,
are denoted by $\hat m_1^{\BF}, \hat m_1^{\SBF}, \hat m_1^{*,\BF}$
and $\hat m_1^{*,\SBF}$. Thus, for simplicity of notation, we use the
same kind of symbol for the updates ($j=1$) and for the inputs of the
backfitting algorithms ($2 \le j \le d$).

We make the following assumptions:
\begin{enumerate}[(A4)]
\item[(A1)] The $d$-dimensional vector $X^i$ has compact support
$I=I_1 \times\cdots\times I_d$ for bounded intervals $I_j=[a_j,b_j]$
and its density $f_X$ is continuous and strictly positive on $I$.
\item[(A2)] There exist constants $C_K,C_S> 0$ such that for all $x_j
\in I_j$, $1\leq j \leq d$, the kernels $K_{j,g}(x_j,\cdot)$ are
positive, bounded by $C_K g^{-1}$, have bounded support $\subset
[x-C_Sg,x+C_Sg]$, and are Lipschitz continuous with Lipschitz constant
bounded by $C_Kg^{-2}$. The weight functions $w_j$ are bounded
functions with $w_j(x_j) \geq0$ for $x_j \in I_j$ and $\int w_j(x_j)
\,dx_j > 0$.
\item[(A3)] The conditional density $f_{\varepsilon|X}(0|x)$ of
$\varepsilon$ given $X=x$ is bounded away from zero and infinity for
$x\in I$. Furthermore, it satisfies the following uniform Lipschitz condition:
\[
|f_{\varepsilon|X}(e|x)-f_{\varepsilon|X}(0|x)| \leq C_1 |e|
\]
for $x\in I$ and for $e$ in a neighborhood of 0 with a constant $C_1>0$
that does not depend on $x$.
\item[(A4)] The bandwidths $h_1, \ldots,h_d$ are of order $n^{-1/5}$.
\end{enumerate}

Assumptions (A1)--(A4) are standard smoothing assumptions. In
particular, (A2) is fulfilled for convolution kernels with an
appropriate boundary correction.

For the properties of the updated estimators, the estimators at the
preceding iteration step need to fulfill certain regularity conditions.
We will proceed with the following assumptions that are stated for some
constants $0<\rho\leq1$, \mbox{$\Delta_1, \Delta_2, \Delta_3 > 0$} and $0
\leq\xi\leq(1+\rho) \Delta_1$.

\begin{enumerate}[(A6)]
\item[(A5)] For $j=2,\ldots,d$, it holds for $l=\BF$ and $l=\SBF$ that
\begin{eqnarray*}
{\sup_{a_j+ C_Sh_j\leq x_j\leq b_j- C_Sh_j}} | \hat
m_j^{l}(x_j)-m_j(x_j) | &=& O_P\bigl(n^{-{(4+4\rho)/(10+15\rho)}
-\Delta_1}\bigr),\\
{\sup_{a_j\leq x_j\leq b_j}} | \hat
m_j^{l}(x_j)-m_j(x_j) | &=& O_P\bigl(n^{- [{(4+4\rho)/(10+15\rho
)}-\Delta_1 ]/2}\bigr).
\end{eqnarray*}
\item[(A6)] There exist random functions $g_2,\ldots,g_d$ with
derivatives that fulfill the Lipschitz condition
\[
|g^\prime_j(x_j)-g^\prime_j(x_j^*)|\leq C |x_j-x_j^*|^\rho n^{\xi}
\]
for $j=2,\ldots,d$ and $x_j,x_j^*\in I_j$. Furthermore, these
functions satisfy
\[
{\sup_{a_j\leq x_j\leq b_j} }| \hat
m_j^{l}(x_j)-g_j(x_j) | = O_P(n^{-2/5-\Delta_2})
\]
for $l=\BF$ and $l=\SBF$.
\item[(A7)] For $j=2,\ldots,d$, it holds for $l=\BF$ and $l=\SBF$ that
\begin{eqnarray*}
{\sup_{a_j+ C_Sh_j\leq x_j\leq b_j- C_Sh_j}} | \hat
m_j^{l}(x_j)-\hat m_j^{*,l}(x_j) | &=& O_P(n^{-2/5-\Delta_3}),\\
{\sup_{a_j\leq x_j\leq b_j} }| \hat
m_j^{l}(x_j)-\hat m_j^{*,l}(x_j) | &=& O_P(n^{-1/5-\Delta_3}).
\end{eqnarray*}
\end{enumerate}

We briefly comment on the assumptions (A5)--(A7). A more detailed
discussion is given after Theorem \ref{TH:2.2}. Assumption (A5) requires
suboptimal rates for the preceding estimators that are plugged in for
the update of the first component. Assumption (A6) states that the
class of possible realizations of the preceding estimators is not too
rich. We assume that the preceding estimators are in a neighborhood of
the class of functions with Lipschitz continuous derivatives. Other
classes could be used but for a Lipschitz class it is relatively easy
to check if a function belongs to it. Note that we do not assume that
the quantile estimator itself has a smooth derivative. In general, such
an assumption does not hold because quantile kernel estimators are not
smooth. Assumption (A7) is very natural. It states that the estimators
that are plugged into the updating equation of the quantile model and
of the mean regression model differ only by second order terms. Without
this assumption, it cannot be expected that the updated estimators
differ also only by second order terms. We will see below that this
assumption is automatically fulfilled if we apply Proposition \ref{TH:2.1} for
an analysis of iterative applications of the backfitting algorithms. In
the assumptions (A5) and (A7), if one replaces the interior region
$[a_j+C_S h_j, b_j-C_S h_j]$ by the whole range $[a_j,b_j]$ and if one
uses boundary corrected kernels, then one can also replace in
Proposition \ref{TH:2.1} the suprema over the interior region by those
over the whole range, and the estimators achieve the rate $n^{-2/5}$ at
the boundary, too.
\begin{proposition} \label{TH:2.1} Under the assumptions \textup{(A1)--(A7)}, it
holds for the updated estimators with $l=\BF$ and with $l=\SBF$ that for
some $\delta>0$
\begin{eqnarray*}
{\sup_{a_1+ C_Sh_1\leq x_1\leq b_1- C_Sh_1}} | \hat
m_1^{l}(x_1)-\hat m_1^{*,l}(x_1) | &=& O_P(n^{-2/5-\delta}),\\
{\sup_{a_1\leq x_1\leq b_1}} | \hat
m_1^{l}(x_1)-\hat m_1^{*,l}(x_1) | &=& O_P(n^{-1/5-\delta}).
\end{eqnarray*}
\end{proposition}

The additional factor $n^{-\delta}$ allows an iterative application of
the proposition. This has an important implication. We recall that the
backfitting algorithms for mean regression have a geometric rate of
convergence. In particular, in the case of smooth backfitting, only
square integrability for the initial estimator is required for the
algorithm to achieve the geometric rate of convergence, see Theorem 1
of \citet{MLN99}. Suppose one chooses square
integrable functions, say $\hat m_2^{\BF,[0]}, \ldots, \hat
m_d^{\BF,[0]}$ as the starting value in the algorithm for the
backfitting quantile estimator and that one runs a cycle of backfitting
iterations (\ref{it:BF}) for $j=1,\ldots,d$. Then we get updates $\hat
m_2^{\BF,[l]}, \ldots,\hat m_d^{\BF,[l]}$ with $l=1$ and after further
cycles with $l >1$. (Note that by construction of the backfitting
estimator we do not need a pilot version of $m_1^{\BF,[0]}$.) Then, one
can think of running the backfitting mean regression algorithm (\ref
{it:TBF}) with the same initial estimators $\hat m_2^{\BF,[0]}, \ldots
,\hat m_d^{\BF,[0]}$ in parallel with the backfitting quantile
regression algorithm~(\ref{it:BF}). This results in updates $\hat m_2^{*,\BF,[l]}, \ldots
,\hat m_d^{*,\BF,[l]}$ for $l\geq1$. In the proof of our next theorem,
we will see that after $l$ cycles of the two parallel iterations, the
difference $\hat m_j^{\BF,[l]}-\hat m_j^{*,\BF,[l]}$ is of order
$O_P(n^{-2/5-\delta})$ in the interior, and of order
$O_P(n^{-1/5-\delta})$ at the boundaries. This holds as long as $l
\leq C_{\mathrm{iter}} \log n$ with $C_{\mathrm{iter}}$ small enough. On the
other hand,\vspace*{-2pt} we will show that $\hat m_j^{*,\BF,[C_{\mathrm{iter}} \log n]}$
is asymptotically equivalent to the limit of the backfitting algorithm
$\hat m_j^{*,\BF,[\infty]}$, if $C_{\mathrm{iter}}$ is large enough. If
the pilot estimators $\hat m_2^{\BF,[0]}, \ldots, \hat m_d^{\BF,[0]}$ are
accurate enough, then the constant $C_{\mathrm{iter}}$ can be chosen such
that both requirements are fulfilled. This will allow us to get the
asymptotic limit distribution of $\hat m_j^{*,\BF,[C_{\mathrm{iter}} \log
n]}$, and thus that of $\hat m_j^{\BF,[C_{\mathrm{iter}} \log n]}$.

Similar findings also hold for the smooth backfitting estimator. We
denote the starting values by $\hat m_2^{\SBF,[0]}, \ldots,\hat
m_d^{\SBF,[0]}$ and the updates by $\hat m_2^{\SBF,[l]}, \ldots,\break\hat
m_d^{\SBF,[l]}$ or $\hat m_2^{*,\SBF,[l]}, \ldots,\hat m_d^{*,\SBF,[l]}$,
respectively.

The following theorem summarizes our discussion. For the theorem, we
need the following additional assumptions:

\begin{enumerate}[(A8)]
\item[(A8)] There exist constants $c_K, C_D> 0$, $C_S^{\prime}\geq0$
such that for $a_j + C_S^{\prime} h_j \leq x_j,u_j\leq b_j-C_S^{\prime
} h_j$ it holds that $K_{j,h_j}(x_j,u_j) = h_j^{-1}
K[h_j^{-1}(x_j-u_j)]$ for a function $K$ with $\int K(v) \,dv = 1$ and
$\int v K(v) \,dv = 0$. For all $x_j,u_j\in I_j$, $1\leq j \leq d$, the
kernels\vspace*{1pt} $K_{j,g}(x_j,u_j)$ have a second derivative w.r.t. $x_j$ that
is bounded by $C_D g^{-3}$ and they fulfill $\int K_{j,g}(x_j,v_j)
\,dv_j \geq c_K$ and $\int K_{j,g}(v_j,u_j) \,dv_j =1$.

\item[(A9)] The function $ f_{X_k|X_j}^w (x_k|x_j) \equiv
f_{X_j,X_k}^w(x_j,x_k)/f_{X_j}^w(x_j)$ has a second\break derivative w.r.t.
$x_j$ that is bounded over $x_j\in I_j$, $x_k\in I_k$, $1\leq j,k \leq
d$, $k \not= j$.
\end{enumerate}

The last condition in (A8) implies that the one-dimensional kernel
density estimators integrate to one and that they are equal to the
corresponding marginalization of higher-dimensional product-kernel
density estimators. This assumption simplifies bias calculation of the
backfitting estimators.
\begin{theorem} \label{TH:2.2} Assume\vspace*{1pt} that \textup{(A1)--(A4), (A8)} and
\textup{(A9)} hold, and that \textup{(A5)} and \textup{(A6)} are satisfied by $\hat m_j^{\BF}=\hat
m_j^{\BF,[0]}$ and $\hat m_j^{\SBF}=\hat m_j^{\SBF,[0]}$ ($j=2,\ldots,d$)
with $\xi,\Delta_2,\Delta_3, {2 \over5} - {1+ \rho\over2 + 3 \rho
}{4 \over5} - \Delta_1 >0$ small enough. Then, we get for $\hat
m_j^{l,\mathrm{iter}}=\hat m_j^{l,[C_{\mathrm{iter}} \log n]}$ with an
appropriate choice of $C_{\mathrm{iter}}=C_{\mathrm{iter},l}$ ($l=\BF$
and $l=\SBF$) that for $a_j < x_j < b_j$
\begin{eqnarray*}
&&
\sqrt{nh_j} [ \hat m_j^{l,\mathrm{iter}} (x_j)
-m_j(x_j) - h_j^2 \beta_j(x_j) ]\\
&&\qquad \to N \biggl(0, {\alpha
(1-\alpha) \over f_{\varepsilon,X_j}(0,x_j)^2}f_{X_j}(x_j) \int
K^2(u) \,du \biggr)
\end{eqnarray*}
in distribution, where $\beta_j(x_j) = \beta_j^*(x_j) - \int\beta
_j^*(u_j)w_j(u_j) \,du_j$, $\beta_j^*(x_j)= h_j^{-2}\times\break m_j^{\prime}
(x_j) \int(u_j-x_j) K_{j,h_j}(x_j,u_j) \,du_j +\mu_{2,K} {1 \over
2} m_j^{\prime\prime}(x_j) +\mu_{2,K} \beta_j^{**}(x_j)$, $\mu
_{2,K}=\break \int v^2 K(v) \,dv$ and $(\beta_1^{**},\ldots,\beta_d^{**})$ is a
tuple of functions that minimizes
\[
\int\Biggl[ \sum_{j=1}^d \biggl(m_j^{\prime}
(x_j) {\partial f_{\varepsilon,
X}(0,x)/\partial x_j \over f_{\varepsilon, X}(0,x) } - \beta
_j^{**}(x_j) \biggr) \Biggr] ^2 f_{\varepsilon, X}(0,x) \,dx.
\]
\end{theorem}

Note that the first term in the definition of $\beta_j^*$ is of order
$n^{1/5}$ at the boundary but vanishes in the interior of $I_j$.
Because of the norming with the weight function $w_j$, the bias
function $\beta_j$ is shifted from $\beta_j^*$ by $\int\beta
_j^*(u_j)w_j(u_j) \,du_j$. One can estimate the bias and the variance
terms because they only require two-dimensional objects if one
calculates them with the backfitting algorithms.

We now come back to discussion of the assumptions (A5)--(A7). Assumption
(A5) allows that the starting estimators have a suboptimal rate. In
particular, it requires that the starting estimators are consistent.
For example, one could use here orthogonal series estimators, smoothing
splines or sieve estimators. In the simulations, we got good results by
using constant functions as starting values, that is, functions that
are not consistent. For backfitting mean regression, it is known that
every starting value works. Because of the nonlinearity of quantile
regression, we do not expect that such a result can be proved for
quantile regression. In our result, we did not specify the required
rate for the pilot estimator. But, if one does this, we conjecture that
one can get the statement of Theorem \ref{TH:2.2} with pilot estimators that
have much slower rates. For such a theorem, one has to prove a
modification of Proposition \ref{TH:2.1} with the following statement: for the
estimators at the preceding stage of the backfitting algorithms, less
accurate error bounds would suffice to get that the difference between
the backfitting estimators $\hat m_1$ and $\hat m_1^*$ at the current
stage of the algorithm is of higher order than the accuracy of the
preceding estimators. This would allow one to weaken the assumptions on
the rate of the starting estimators.

Assumption (A7) is not required for Theorem \ref{TH:2.2}. This is because
running the iterative algorithms
(\ref{it:TBF}) and (\ref{it:TSBF}) is only imaginary and in the proof
we choose to use the same starting values as in the real iterative
algorithms (\ref{it:BF}) and (\ref{it:SBF}), respectively. Thus, (A7)
is automatically satisfied at the beginning of the iterations.
Proposition \ref{TH:2.1} tells us that the updated estimators also fulfill (A7).
This holds with the same rate but with multiplicative factors. For this
reason, after $L$ backfitting cycles the difference between the mean
regression\vspace*{1pt} and the quantile regression estimators is not of order $(C
\times L) n^{-2/5-\delta}$, but of order $C^L n^{-2/5-\delta}$, for
some $\delta> 0$, $C > 1$. For a number of iterations, $C_{\mathrm{iter}}
\log n$ such that $C_{\mathrm{iter}} \log C < \delta$ this is of order
$o(n^{-2/5})$.

Compared with the results for mean regression backfitting estimators,
our results for quantile estimation are weaker in two aspects. First,
we need initial estimators that are consistent, whereas in mean
regression one can start with arbitrary initial values. This
restriction comes from the nonlinearity of the quantile functional.
Second, we put restrictions on the number of iteration steps. It must
be of logarithmic order with a factor that is not too small and not too
large. When letting run the two parallel backfitting procedures for
mean and quantile regression, we were not able to control in the proof
the difference between the two outcomes if the number of iterations is
too large. We conjecture that both restrictions are necessary only for
technical reasons in our approach for the proof. In our simulation, we
started with nonconsistent pilot estimators and we let the algorithms
run until the outcomes were stabilized. According to our experience in
the simulation, there seemed practically no advantage in limiting the
number of iterations and there was also no problem when starting the
algorithm with initial estimators that were far away from the
corresponding underlying regression functions.

A natural extension of our results is to study local polynomial
quantile estimators. This can be done along the lines of this paper by
putting smoothness restrictions also on the higher order terms of the
local polynomial fit. This can be done relatively easily for local
polynomial smooth backfitting. For local polynomial ordinary
backfitting, it would require also essentially new theoretical results
for mean regression. We do not follow this line in this paper.

%s3 ###
\section{Numerical implementation}
\label{sec:num}
In practical implementations of the smooth backfitting
method, one may approximate the integral at (\ref{it:SBF}) by Monte
Carlo integration. This can be done in several ways. In one version,
one generates $(U_2^j, \ldots, U_d^j)$ for $1 \le j \le M$ from a
$(d-1)$-variate uniform distribution on
$I_2 \times\cdots\times I_d$. Then an approximation of $\hat
m^{\SBF}_1(x_1)$ may be obtained by
\begin{eqnarray*}
\hat m^{\SBF}_1(x_1) &\approx& \mathop{\arg\min}_{\theta\in\Theta
}\sum
_{i=1}^{n}\sum_{j=1}^M \tau_{\alpha} \bigl(Y_{i}-\theta-\hat
m^{\SBF}_0-\hat m^{\SBF}_2(U_2^j)-\cdots- \hat m^{\SBF}_d(U_d^j)
\bigr)\nonumber\\
&&\hspace*{64pt}{} \times K_{1,h_1}(x_1,X_1^{i})
K_{2,h_2}(U_2^j,X_2^{i}) \cdots
K_{d,h_d}(U_d^j,X_d^{i}).
\end{eqnarray*}
In practical implementation, the values $U_k^j$ can be chosen from a
finite grid of equidistant points.
Then the algorithm has to update the function values of the additive
components on this grid.

In another version, one generates independent $U_{\ell,i,j}$ for $\ell
=2,\ldots,d$, $i=1,\ldots,n$, $j=1,\ldots,J$, where $U_{\ell,i,j}$
has density $ K_{\ell,h_l}(\cdot,X_\ell^{i})$. Again, in practical
implementation, the values of these random variables can be chosen from
a finite grid of equidistant points. Then the smooth backfitting
estimator at $x_1$ is calculated by
\begin{eqnarray*}
\hat m^{\SBF}_1(x_1) &\approx& \mathop{\arg\min}_{\theta\in\Theta
} \sum
_{i=1}^{n}\sum_{j=1}^J \tau_{\alpha} \bigl(Y_{i}-\theta-\hat
m^{\SBF}_0 -\hat m^{\SBF}_2(U_{2,i,j})\\
&&\hspace*{79.8pt}\hspace*{40.5pt}{} -\cdots- \hat m^{\SBF}_d(U_{d,i,j})
\bigr)K_{1,h_1}(x_1,X_1^i).
\end{eqnarray*}
This means that the smooth backfitting estimator can be calculated by
an algorithm that is designed for the ordinary backfitting with sample
$(Y_i,X_1^i,\break U_{2,i,j},\ldots, U_{d,i,j})$ for $i=1,\ldots,n$ and
$j=1,\ldots,J$. In this case, the speed of the algorithm for the
smooth backfitting behaves as that for the ordinary backfitting with
sample size~$Jn$.

In the last algorithm, the values $U_{\ell,i,j}$ could be replaced by
deterministic choices such that for fixed $i$ and $\ell$ the
probability density $ K_{\ell,h_\ell}(\cdot,X_\ell^{i})$ put equal
mass between neighbored points of $U_{\ell,i,j}$, that is,
\[
\int_{-\infty}^{U_{\ell,i,j}}
K_{\ell,h_\ell}(x_\ell,X_\ell^{i}) \,dx_\ell= j/(J+1),\qquad
j=1,\ldots,J.
\]
Suppose that $K_{\ell,h_\ell}(\cdot, z)$ is symmetric about $z$.
Then the algorithm calculates the ordinary backfitting estimates when
$J=1$, since in that case $U_{\ell,i,1} = X_\ell^{i}$. It also
approximates the smooth backfitting estimates as $J\to\infty$. Thus,
there exists a broad band of compromises between the ordinary
backfitting and the smooth backfitting for intermediate choices of $J$.

%%%%%%%%% Section 4
%s4 ###
\section{Simulation study}
\label{sec:sim}

In this section, we illustrate the asymptotic equivalence asserted in
Proposition \ref{TH:2.1}. We compared the numerical properties of the ordinary
backfitting (BF) and the smooth backfitting (SBF) estimators defined at
(\ref{it:BF}) and (\ref{it:SBF}) with their theoretical mean
regression versions defined at (\ref{it:TBF}) and (\ref{it:TSBF}),
respectively.

In the simulation, we considered the following model:
\[
Y^i = f_1(X_1^i) + f_2(X_2^i) + f_3(X_3^i)
+\{\sigma_1(X_1^i) + \sigma_2(X_2^i)+ \sigma_3(X_3^i)\} U^i,
\]
where $U^i$ are i.i.d. $N(0,1)$, $f_1(x_1) = x_1^3$, $f_2(x_2) =\sin
(\pi x_2)$, $f_3(x_3) = 2\times\break \exp(-16x_3^2)$, $\sigma_1(x_1)=\cos(x_1)$,
$\sigma_2(x_2) = \exp(x_2)$ and $\sigma_3(x_3)=\exp(x_3)$. With
this model, the centered version of the $j$th additive component of the
$\alpha$-quantile function equals
\[
m_j (x_j; \alpha) = c_j + f_j(x_j) + \sigma_j(x_j)
\Phi^{-1}(\alpha),
\]
where $\Phi^{-1}(\alpha)$ is the $\alpha$-quantile of the standard
normal distribution and $c_j$ is the constant that makes $E m_j (X_j^1;
\alpha) =0$. We considered two different cases for the distribution of
$X^i$. One was the case where the components of $X^i$ were independent.
In this case, $X^i$ were generated from $N_3(0, J)$ truncated outside
$[-1,1]^3$, where $J$ denotes the identity matrix of dimension $d=3$.
This means the density of $X^i$ was $f_X(x) = \varphi(x)I(x \in
[-1,1]^3)/\int_{[-1,1]^3} \varphi(z) \,dz$, where $\varphi$ denotes
the density function of $N_3(0, J)$. The second was the case where the
components of $X^i$ were correlated. In this case, $X^i \sim N_3(0, V)$
truncated outside $[-1,1]^3$, where $V \equiv(v_{ij})$ has $v_{ii}=1$
and $v_{ij}=0.9$ for $i \neq j$. Because of the truncation, the actual
correlation equals $0.644$. The sample sizes were $n=200$ and $n=500$.
These relatively large sample sizes were considered to let the
asymptotic results in Section~\ref{sec:2} be well in effect.

Implementation of the ordinary and smooth backfitting methods requires
optimization involving the nonsmooth function $\tau_\alpha$. For
this, we used R function \texttt{rq()} in the library \texttt{quantreg}. For
the smooth backfitting, we discretized the integrals on a fine grid in
$[-1,1]^3$. We used
%
%e4.1 ###
%
\begin{equation} \label{kernelchoice}
K_{j,g}(x,u) = \biggl[\int K \biggl(\frac{x-u}{g} \biggr) \,dx
\biggr]^{-1} K \biggl(\frac{x-u}{g} \biggr),
\end{equation}
where $K$ is Epanechinikov kernel given by $K(u) =
(3/4)(1-u^2)I_{[-1,1]}(u)$. For the bandwidths, we took $h_1 = h_2 =h_3
= h$ for simplicity. Normalization was done in each iteration so that
$\int\hat m_j (x_j)\hat f_{X_j}(x_j) \,dx_j =0$. Note that we used
estimates of $f_{X_j}$ in the normalization, instead of fixed weight
functions which we considered in our theoretical development for
simplicity. Using a different weight function changes the estimator
only by an additive constant. To get the density estimates~$\hat
{f}_{X_j}$, we used the same kernel $K$ and the bandwidth $h$ that we
employed for quantile estimation. We chose the initial estimates in the
iterative algorithms (\ref{it:BF}),
(\ref{it:SBF}), (\ref{it:TBF}) and (\ref{it:TSBF}) to be zero. It
was found that the algorithms converged with this initial choice in all cases.

%t1 ###
%
\begin{table}
\caption{Mean integrated squared errors of the
estimators}\label{table1}
\begin{tabular*}{\tablewidth}{@{\extracolsep{\fill}}l c c c c c@{}}
\hline
\textbf{Sample} & \textbf{Distribution} &&&&\\
\textbf{size} & \textbf{of} $\bolds X$ & \textbf{Method}
& $\bolds{\alpha=0.2}$ & $\bolds{\alpha=0.5}$ & $\bolds{\alpha=0.8}$
\\
\hline
$n=200$&Uncorrel.& BF & 0.09345& 0.07457& 0.08770\\
& & BF$^*$& 0.09585& 0.07512& 0.08208\\
& & SBF& 0.08818& 0.07039& 0.08209\\
& & SBF$^*$& 0.09436& 0.07455& 0.07937\\
& Correl.& BF& 0.09043& 0.07165& 0.08382\\
& & BF$^*$& 0.09864& 0.07539& 0.08276\\
& & SBF& 0.08555& 0.06712& 0.07937\\
& & SBF$^*$& 0.09136& 0.07140& 0.08412\\
[4pt]
$n=500$&Uncorrel.& BF& 0.05240& 0.04020& 0.04881\\
& & BF$^*$& 0.04959& 0.04121& 0.04729\\
& & SBF& 0.04905& 0.03827& 0.04557\\
& & SBF$^*$& 0.05045& 0.04178& 0.04896\\
& Correl.& BF& 0.05463& 0.04182& 0.05094\\
& & BF$^*$& 0.05137& 0.04305& 0.05312\\
& & SBF& 0.05186& 0.03983& 0.04743\\
& & SBF$^*$& 0.05496& 0.04221& 0.05296\\
\hline
\end{tabular*}
\legend{Note: BF$^*$ denotes the theoretical mean
regression ordinary backfitting estimator, and SBF$^*$ denotes the
theoretical mean regression smooth backfitting estimator.}
\end{table}

Table \ref{table1} show Monte Carlo estimates, based on 200
pseudo-samples, of the mean integrated squared errors,
\[
\mathrm{MISE} = E \int\{\bar m_1(x_1) + \bar m_2(x_2) +
\bar
m_3(x_3) -m_1(x_1)- m_2(x_2) - m_3(x_3) \}^2 f_X(x) \,dx,
\]
where $f_X$ is the density function of $X^i$, and $\bar m_j$
represents $\hat m_j^{\BF}$, $\hat m_j^{\SBF}$, $\hat m_j^{*,\BF}$ or
$\hat m_j^{*,\SBF}$. For each\vspace*{1pt} estimator, its MISE was estimated by
$\overline{\mathrm{ISE}} = \sum_{r=1}^{200} \mathrm{ISE}_r/200$, where
$\mathrm{ISE}_r$ is the value of the integrated squared error
\[
\int\{\bar m_1(x_1) + \bar m_2(x_2) + \bar
m_3(x_3) -m_1(x_1)- m_2(x_2) - m_3(x_3) \}^2 f_X(x) \,dx
\]
for the $r$th sample. We computed the estimates of the additive
regression function with bandwidths on a grid in $[0.1, 1.5]$. The
values for $\hat m^{\BF}$ and $\hat m^{*,\BF}$ reported in the table are
for the bandwidths that gave optimal performance of $\hat m^{\BF}$, and
likewise those for $\hat m^{\SBF}$ and $\hat m^{*,\SBF}$ are for the
bandwidths that gave optimal performance of $\hat m^{\SBF}$. In most
cases, the estimated MISE was minimized around $h=0.5$ when $n=200$,
and around $h=0.4$ when $n=500$. This is roughly consistent with the
theory that the size of the optimal bandwidth equals $n^{-1/5}$ for
univariate smoothing, according to which the ratio of the optimal
bandwidths for $n=200$ and $n=500$ equals $(500/200)^{1/5} \approx1.20$.

To compare $\hat m^{\BF}$ and $\hat m^{\SBF}$ with their theoretical mean
regression counterparts $\hat m^{*,\BF}$ and $\hat m^{*,\SBF}$, we find
that the two corresponding MISE values are very close, and that in most
cases the differences get smaller as $n$ increases. This supports our
theory presented in Section \ref{sec:2}. In the table, we also find that the size
of the estimated MISE for $n=500$ is nearly half of the corresponding
value for $n=200$. This supports the fact that the ordinary and smooth
backfitting estimators enjoy the univariate rate of convergence
$n^{-4/5}$ in MISE, since $(500/200)^{4/5} \approx2.08$.

According to Table \ref{table1}, the MISE values of the estimators at
$\alpha=0.5$ are always smaller than those at $\alpha=0.2$ and
$\alpha=0.8$. Note that, in Theorem \ref{TH:2.2}, $f_{X_j}^w(x_j)$ is
nothing else than the joint density of $(\varepsilon, X_j)$ at the
point $(0,x_j)$. Under our simulation model, the conditional density
can be expressed as
\begin{eqnarray*}
f_{X_j}^w(x_j) &=& \int\frac{1}{\sigma_1(x_1)+\sigma_2(x_2)+\sigma
_3(x_3)}\phi\biggl( \frac{\Phi^{-1}(\alpha)}{\sigma_1(x_1)+\sigma
_2(x_2)+\sigma_3(x_3)} \biggr)\\
&&\hspace*{8.8pt}{}\times f_{X}(x) \,dx_{-j}
\end{eqnarray*}
for $j=1, 2$ and $3$, where $\phi$ denotes the density of the standard
normal distribution. According to Theorem \ref{TH:2.2}, this implies
that the theoretical value of the integrated variance increases as
$\alpha$ gets away from $0.5$. This explains why we have larger MISE
values for $\alpha$ away from $0.5$. Similar numerical evidences were
also observed by \citet{YJ98} and \citet{LLP06}.

%f1 ###
%
\begin{figure}

\includegraphics{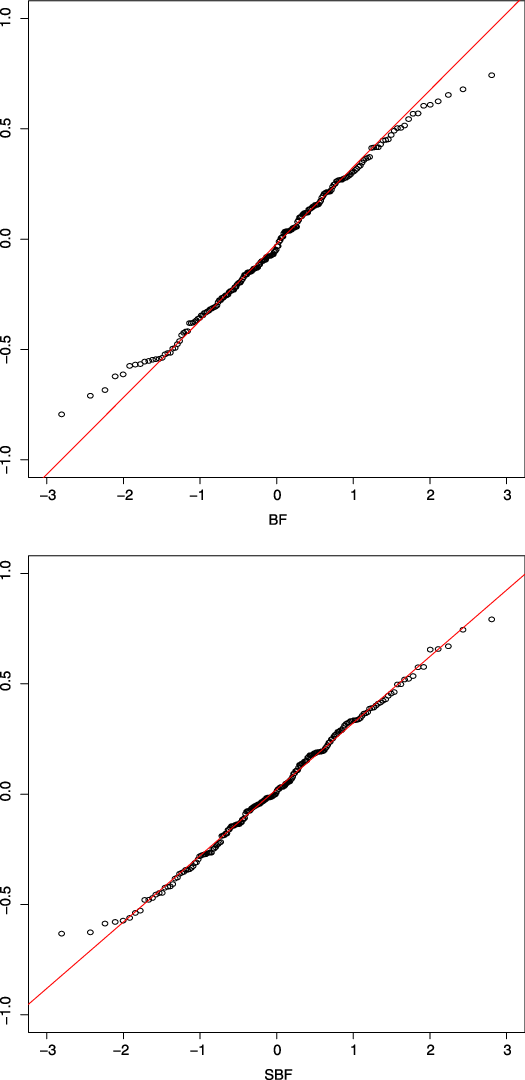} %[scale=0.96]

\caption{Normal Q--Q plots for $\hat m_2^{\BF}(x)$ and
$\hat m_2^{\SBF}(x)$ based on 200 values computed from pseudo-samples in
the case where $x=0$, $\alpha=0.5$, $n=200$ and the components of $X^i$
were correlated. The theoretical quantiles are on the horizontal axis
and the sample quantiles are on the vertical axis.}
\label{fig1}
\end{figure}

Figure \ref{fig1} illustrates the asymptotic normality of $\hat
m_j^{\BF}$ and $\hat m_j^{\SBF}$. It depicts the normal Q--Q plots of the
200 values of $\hat m_2^{\BF}(x)$ and $\hat m_2^{\SBF}(x)$ at $x=0$ when
$\alpha=0.5$ and $n=200$. The figure is for the case where the
components of $X^i$ are correlated. Although it exhibits slight
departures from normality at tails, the figure suggests that the
distributions of the estimators get close to normal even for moderate
sample sizes. We obtained other Q--Q plots that corresponded to other
components $j$, other points $x$ or other quantile levels $\alpha$, and
also repeated them in other simulation models. They looked not much
different from the case we report here.

%f2 ###
%
\begin{figure}

\includegraphics{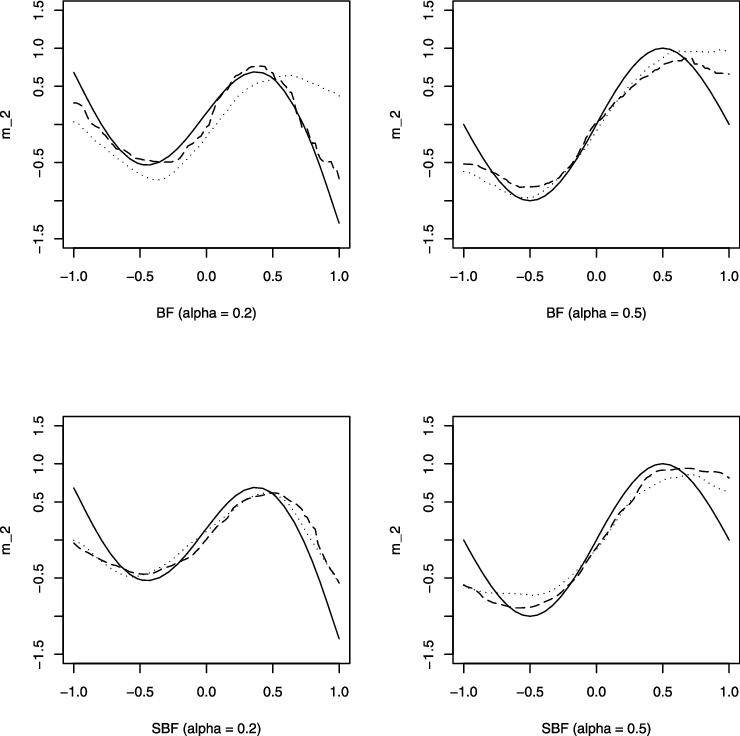}

\caption{Estimates of a component function computed from a
sample that gave the median performance in terms of the integrated
squared\vspace*{-2pt} error of $\hat{m}_j^{\BF}$ or $\hat{m}_j^{\SBF}$,
when $n=500$ and the covariates were correlated.
Long-dashed\vspace*{-2pt} and dotted curves in the top two panels are
$\hat{m}_j^{\BF}$ and $\hat{m}_j^{*,\BF}$, respectively, and
those in the bottom two panels are $\hat{m}_j^{\SBF}$ and
$\hat{m}_j^{*,\SBF}$. Left two panels are for the case
$\alpha=0.2$ and the right are for $\alpha=0.5$. Solid curves
represent the true component functions.}
\label{fig2}
\end{figure}

Figure \ref{fig2} illustrates\vspace*{1pt} how the four curve estimates $\hat
m_j^{\BF}$, $\hat m_j^{*,\BF}$, $\hat m_j^{\SBF}$ and $\hat m_j^{*,\SBF}$
computed from a single typical sample look like. In the top two panels,
the long-dashed and dotted curves, respectively, represent $\hat
m_j^{\BF}$ and $\hat m_j^{*,\BF}$ computed from a sample for which the
value of the integrated squared error
\[
\int\{\hat m_j^{\BF}(x_j) - m_j(x_j) \}^2 \,dx_j
\]
was the median of those values obtained from the 200 pseudo-samples.
Similarly, the bottom two panels depict $\hat m_j^{\SBF}$ and $\hat
m_j^{*,\SBF}$ computed from a sample that gave the median performance in
terms of the integrated squared error
\[
\int\{\hat m_j^{\SBF}(x_j) - m_j(x_j) \}^2 \,dx_j.
\]
In the figure the solid curves represent the true functions. In
comparison of the pairs, $m_j^{\BF}$ versus $\hat m_j^{*,\BF}$ and $\hat
m_j^{\SBF}$ versus $\hat m_j^{*,\SBF}$, we find that the two
corresponding curves move together relatively closer than with the true
function, although there are some places where they are more distant in
the case of the backfitting estimator for $\alpha=0.2$ (top left
panel). The figure is for the estimates of the second component
function when $n=500$ and the components of $X^i$ were correlated.
Those for other cases gave similar lesson, so that are not included here.

%t2 ###
%
\begin{table}
\caption{Differences in mean integrated squared errors
of BF and SBF estimators}\label{table2}
\begin{tabular*}{\tablewidth}{@{\extracolsep{\fill}}l c c c c@{}}
\hline
\textbf{Sample} & \textbf{Distribution} &&&\\
\textbf{size} & \textbf{of} $\bolds X$ & $\bolds{\alpha=0.2}$
& $\bolds{\alpha=0.5}$ & $\bolds{\alpha=0.8}$\\
\hline
$n=200$&Uncorrel.& 0.00527& 0.00418& 0.00561\\
& & (0.00099)& (0.00063)& (0.00087)\\
& Correl.& 0.00488& 0.00453& 0.00445\\
& & (0.00098)& (0.00068)& (0.00096)\\[4pt]
$n=500$&Uncorrel.& 0.00335& 0.00193& 0.00324\\
& & (0.00045)& (0.00028)& (0.00038)\\
& Correl.& 0.00277& 0.00199& 0.00351\\
& & (0.00042)& (0.00034)& (0.00042)\\
\hline
\end{tabular*}
\legend{Note: the numbers are averages of $(\mbox{ISE of
}\hat m^{\BF}) - (\mbox{ISE of }\hat m^{\SBF})$ over 200
pseudo-samples, and their standard errors are given in the
parentheses.}
\end{table}

One may be also interested in comparing the two backfitting quantile
estimators $\hat m^{\BF}$ and $\hat m^{\SBF}$ in terms of MISE. For this,
we computed the standard errors of the differences between the
estimated values of MISE of the respective estimators. In Table \ref
{table2}, we
provide the average differences $\overline{\mathrm{DIFF}}$ and their
standard errors calculated by the formula
\[
\mathrm{S.E.} = \sqrt{\sum_{r=1}^{200} (\mathrm{DIFF}_r -
\overline{\mathrm{DIFF}})^2/(199\times200)},
\]
where $\overline{\mathrm{DIFF}}$ denotes the average of
$\mathrm{DIFF}_r$ over
200 pseudo-samples, and
\[
\mathrm{DIFF}_r = (\mbox{ISE of } \hat m^{\BF} \mbox{ for the }
r\mbox{th sample}) - (\mbox{ISE of } \hat m^{\SBF} \mbox{ for the
} r\mbox{th sample}).
\]

Comparing the two backfitting quantile estimators, we find that the
smooth backfitting estimators have smaller values of the estimated MISE
in all cases than the ordinary backfitting estimators. In particular,
all the differences are statistically significant, exceeding two
standard errors. Although not reported in the paper, we also compared
the two backfitting quantile estimators with their oracle versions. An
oracle estimator of an additive component is the one obtained by using
true functions for the other components. We found that in all cases the
two backfitting quantile estimators had similar performance as their
oracle versions.

%%%%%%%%% Section 4
%s5 ###
\section{Proofs}\label{sec5}

%s5.1 ###
\subsection[Proof of Proposition 2.1]{Proof of Proposition \protect\ref{TH:2.1}}
\label{sec:4.1}

We only give the proof for the ordinary backfitting etimator. The proof
will be given for $a_1+C_S n ^{-1/5} \leq x_1 \leq b_1-C_S n ^{-1/5}$.
The proofs for the smooth backfitting estimator and for boundary points
follow by similar arguments. For simplicity of notation, we also assume
that $d=2$.

The basic asymptotic argument for a treatment of parametric and
nonparametric quantile estimators is a Bahadur expansion. It states
that the quantile estimator is asymptotically equivalent to a linear
statistic, that is, to a sum of independent variables. This expansion
would directly carry over to our case if the pilot functions (input) of
the backfitting algorithms would be nonrandom. Because this is not the
case, we have to generalize the Bahadur approach. We have to show that
the Bahadur expansion holds uniformly over a class of pilot functions.
Furthermore, we have to verify that the pilot estimators lie in this
function class with probability tending to one. The latter is
guaranteed by the assumptions (A5) and (A6). The uniform expansion is the main
step of our proof.

Define
\begin{eqnarray*}
&&
V_i(\theta,\mu_2,x_1) \\
&&\qquad= K_{1,h_1}(x_1,X_1 ^i) \bigl[\tau_\alpha
\bigl(Y^i -\theta- \mu_2(X_2 ^i) \bigr) - \tau_\alpha
\bigl(\varepsilon^i + m_1(X^i_1) - m_1(x_1) \bigr)\\
&&\qquad\quad\hspace*{62pt}{} - \bigl(\theta-m_1(x_1)+\mu_2(X_2^i)-
m_2(X_2^i) \bigr)\\
&&\quad\qquad\hspace*{118.4pt}{}\times\bigl(I\bigl(\varepsilon^i+m_1(X^i_1)-m_1(x_1) < 0
\bigr)-\alpha\bigr) \bigr].
\end{eqnarray*}
Let $J_1 \equiv J_1(x_1)$ and $J_2 \equiv J_2(x_1)$ be index
sets defined by
\begin{eqnarray*}
J_1&=& \{ i\dvtx |X_1^i-x_1|\leq C h_1, a_2 + C_S n^{-1/5}
\leq X_2^i \leq b_2 - C_S n^{-1/5}\},\\
J_2&=& \{ i\dvtx |X_1^i-x_1|\leq C h_1, a_2 \leq X_2^i < a_2 + C_S
n^{-1/5} \mbox{ or } b_2 - C_S n^{-1/5} < X_2^i \leq b_2\}.
\end{eqnarray*}
Put
\begin{eqnarray*}
&&D(\theta,\mu_2,x_1) \\
&&\qquad = \sum_{i=1}^n [V_i(\theta,\mu_2,x_1)-
E^{\mathcal{X}} V_i(\theta,\mu_2,x_1) ] \\
&&\qquad = \sum_{i\in J_1} [V_i(\theta,\mu_2,x_1)-
E^{\mathcal{X}}
V_i(\theta,\mu_2,x_1) ]\\
&&\qquad\quad{} + \sum_{i\in J_2}
[V_i(\theta,\mu_2,x_1)- E^{\mathcal{X}} V_i(\theta,\mu_2,x_1)]\\
&&\qquad \equiv D_1(\theta,\mu_2,x_1)+D_2(\theta,\mu_2,x_1),
\end{eqnarray*}
where $E^{\mathcal{X}}$ is the conditional expectation given $\mathcal
{X} =\{X^1, \ldots, X^n\}$.
Let $M_1$ and $M_2$ denote the numbers of elements of $J_1$ and $J_2$,
respectively. These are random variables. Since $h_1$ is of order
$n^{-1/5}$ and the density $f_X$ is strictly positive on its
support,\vspace*{1pt}
$M_1$ is of order $n\times n^{-1/5}=n^{4/5}$ and $M_2$ is of order
$n\times n^{-1/5}\times n^{-1/5}=n^{3/5}$. Thus, there exist constants
$C_1>0$ and $C_2 > 0$ such that $C_1 n^{4/5} \leq M_1 \leq2 C_1
n^{4/5}$ and $C_2 n^{3/5} \leq M_2 \leq2 C_2 n^{3/5}$ with
probability tending to one.

For a fixed constant $D>0$, we now introduce the class $\mathcal{M}_n
$ of all tuples of a parameter $\theta\in\Theta$ and a function $g$
that fulfills
\[
{\sup_{a_2+ C_Sn^{-1/5}\leq x_2\leq b_2- C_Sn^{-1/5}}} |
g(x_2)-m_2(x_2) | \leq D n^{-{(1+\rho)/(2+3\rho)}{4 /
5}-\Delta_1}
\]
and whose derivative fulfills a Lipschitz condition of order $\rho$
with Lipschitz constant $C$ as in (A6).

For $j\geq0$, let $\mathcal{M}_n(2^{-j})$ denote a grid of points in
$\mathcal{M}_n $ such that for every $(\theta,g) \in\mathcal{M}_n $
there exists $(\theta^*,g^*) \in\mathcal{M}_n(2^{-j}) $ with
$|\theta^*- \theta|\leq2^{-j}$ and $\|g^*- g\|_\infty\leq2^{-j}$.
Let $N_j$ denote the number of points in the grid $\mathcal
{M}_n(2^{-j})$. Note that $N_j= O \{\exp(2^{j/(1+\rho)}n^{\xi
/(1+\rho)}) \}$.

We apply the Bernstein inequality. For a sum of $r$ independent random
variables $V_i$ that are absolutely bounded by a constant $\kappa$ and
have finite variance bounded by $\sigma^2$, this inequality states that
\begin{eqnarray*}
P \Biggl( \Biggl| r^{-1/2} \sum_{i=1}^r (V_i - E V_i) \Biggr| \geq
a \Biggr) &\leq& 2 \exp\biggl( - {a^2 \over2 a \kappa
r^{-1/2} + 2 \sigma^2}
\biggr)\\
&\leq& 2 \exp\biggl( - {a \over4 \kappa r^{-1/2}} \biggr) + 2
\exp\biggl( - {a ^2 \over4\sigma^2} \biggr).
\end{eqnarray*}
We apply this inequality with a chaining argument for $D_1(\theta,\mu
,x_1)$ and $D_2(\theta$,\break $\mu,x_1)$. In doing this, we take $r=M_1$ (or
$r=M_2$, resp.) and $P=P^{\mathcal{X}}$ where $P^{\mathcal
{X}}$ is the conditional distribution given $\mathcal{X} =\{X^1,
\ldots, X^n\}$. Let $J_n$ be chosen so that $2^{-J_n} \leq
n^{-2/5-\delta} \leq2^{-J_n+1}$ with $\delta>0$ small enough, see
below. Define $\gamma= 4(1+\rho)/[5(2+3\rho)]$ and $I_n = \{j\dvtx j\leq
J_n, D n^{-\gamma-\Delta_1} \geq2^{-j}\}$. Furthermore,
for $(\theta,\mu)\in\mathcal{M}_n(2^{-J_n})$ choose $(\theta^j,\mu
^j) \in\mathcal{M}_n (2^{-j})$ with
$|\theta^j-\theta|\leq2^{-j}$ and $\|\mu^j-\mu\|_\infty\leq
2^{-j}$. For $j=J_n$, we choose $(\theta^j,\mu^j)=(\theta,\mu)$. We
do not indicate the dependence of $(\theta^j,\mu^j)$ on $(\theta,\mu
)$ in the notation. For $j\leq j_n= \min I_n$, the grid $\mathcal
{M}_n(2^{-j})$ can be chosen so that it contains only one value of $\mu
$. We assume that this value is equal to $\mu^0=m_2$. Furthermore, we
choose $\theta^0=m_1(x_1)$ and we assume w.l.o.g. that the diameter
of $\Theta$ is less than one. For $j=0$, the grid $\mathcal
{M}_n(2^{-j})$ contains only one value which we choose to be $(\theta
^0,\mu^0)$. Then
\begin{eqnarray*}
&&P \Bigl( {\sup_{(\theta,\mu)\in\mathcal{M}_n(2^{-J_n})}}
|D_1(\theta,\mu,x_1)| > n^{-4/5- 2 \delta} \big| \mathcal{X}
\Bigr)\\
&&\qquad \leq P \biggl( {\sup_{(\theta,\mu)\in
\mathcal{M}_n(2^{-J_n})}} \biggl| D_1(\theta^0,\mu^0,x_1)\\
&&\qquad\quad\hspace*{79.3pt}{} +\sum_{1
\leq
j < j_n} D_1(\theta^j,\mu^0,x_1)-D_1(\theta^{j-1},\mu^{0},x_1) \\
&&\qquad\quad\hspace*{79.3pt}{} + \sum_{j_n \leq j \leq J_n}
D_1(\theta^j,\mu^j,x_1)-D_1(\theta^{j-1},\mu^{j-1},x_1) \biggr|\\
&&\qquad\quad\hspace*{233.3pt}> n^{-4/5- 2 \delta} \big| \mathcal{X} \biggr).
\end{eqnarray*}
Let $s_j$ be positive numbers (depending on $n$) such that $\sum
_{1\leq j\leq J_n} s_j \leq1/2$. Then the right-hand side of the above
inequality is bounded by
%
%e5.1 ###
%
\begin{eqnarray}\label{bound1}\qquad
&&P \bigl( |D_1(\theta^0,\mu^0,x_1) | > 2^{-1}
n^{-4/5- 2 \delta} | \mathcal{X} \bigr) \nonumber\\
&&\qquad{} + \sum_{1 \leq j < j_n} 2^{2j} \sup_{*} P \bigl( |D_1(\theta^j,\mu^0,x_1)\nonumber\\
&&\hspace*{113pt}{}-D_1(\theta
^{j-1},\mu^0,x_1) | > s_j n^{-4/5- 2 \delta} | \mathcal
{X} \bigr) \\
&&\qquad{} + \sum_{j_n \leq j \leq J_n} N_j N_{j-1} \sup_{**}
P \bigl( |D_1(\theta^j,\mu^j,x_1)\nonumber\\
&&\hspace*{138pt}{}-D_1(\theta^{j-1},\mu
^{j-1},x_1) | > s_j n^{-4/5- 2 \delta} |
\mathcal{X} \bigr), \nonumber
\end{eqnarray}
where $\sup_{*}$ and $\sup_{**}$ runs over all $(\theta^{j},\mu
^{j})\in\mathcal{M}_n(2^{-j})$ and
$(\theta^{j-1},\mu^{j-1})\in\mathcal{M}_n(2^{-j+1})$ with $|\theta
^j-\theta^{j-1}|\leq2^{-j+1}$ and $\|\mu^j-\mu^{j-1}\|_\infty\leq
2^{-j+1}$.

Using the Bernstein inequality with $\kappa= O(2^{-j} h_1^{-1})$,
$\sigma^2= 2^{-2j} O(n^{-\gamma-\Delta_1}\times\break h_1^{-2})$ and $a
=M_1^{-1/2} n s_j n^{-4/5- 2 \delta} c$, the last sum in (\ref
{bound1}) can be
bounded by
%
%e5.2 ###
%
\begin{eqnarray}\label{bound2}\qquad
&&\sum_{j_n \leq j \leq J_n} \bigl[\exp\bigl(d_1 2^{ j/(1+\rho)} n^{ \xi
/(1+\rho)}- d_2
s_j n
n^{-4/5- 2 \delta} 2^{j}h_1\bigr) \nonumber\\[-8pt]\\[-8pt]
&&\hspace*{35.4pt}{} + \exp\bigl(d_1 2^{j/(1+\rho)} n^{ \xi/(1+\rho)}-
d_2 s^2_j M_1^{-1} n^2
n^{-8/5- 4 \delta} 2^{2j}n^{\gamma+\Delta_1} h^2_1\bigr) \bigr] \nonumber
\end{eqnarray}
for some constants $d_1,d_2 >0$. Choosing $s_j = (d_3 \log n)^{-1}$
with $d_3$ large enough, the sum at (\ref{bound2}) can be bounded
further by
\[
\exp(-d_4 n^{d_5} ) + \exp(-d_6
M_1^{-1} n^{4/5+d_7} ),
\]
where $d_4,\ldots,d_7>0$ are some constants. Here, we used that
$\delta> 0$ is small enough.
Using similar arguments for the first two terms in (\ref{bound1}), one
can bound the sum of all three terms in (\ref{bound1}) by
\[
\exp(-d_8 n^{d_9} ) ,
\]
where $d_8,d_9>0$ are some constants. This exponential bound entails
that for $\delta>0$ small enough
%
%e5.3 ###
%
\begin{eqnarray} \label{uni1}
&&{\mathop{\sup_{(\theta,\mu_2)\in\mathcal{M}_n(2^{-J_n})}}_{x_1 \in
I_1}}
|n^{-1}D_1(\theta,\mu_2,x_1)|\nonumber\\
&&\qquad = {\mathop{\sup_{(\theta,\mu_2)\in\mathcal{M}_n(2^{-J_n})}}_{x_1 \in
I_1}}
\biggl|n^{-1}\sum_{i\in J_1} \{V_i(\theta,\mu_2,x_1)-
E^{\mathcal{X}} V_i(\theta,\mu_2,x_1) \} \biggr|\\
&&\qquad = O_P(n^{-4/5-\delta}). \nonumber
\end{eqnarray}
Similarly, it can be shown that
%
%e5.4 ###
%
\begin{eqnarray} \label{uni2}
&&{\mathop{\sup_{(\theta,\mu_2)\in\mathcal{M}_n(2^{-J_n})}}_{x_1 \in
I_1}}|n^{-1}D_2(\theta,\mu_2,x_1)|\nonumber\\
&&\qquad = \mathop{\sup_{(\theta,\mu_2)\in\mathcal{M}_n(2^{-J_n})}}_
{x_1 \in I_1}
\biggl|n^{-1}\sum_{i\in J_2} \{V_i(\theta,\mu_2,x_1)-
E^{\mathcal{X}}
V_i(\theta,\mu_2,x_1) \} \biggr| \\
&&\qquad = O_P(n^{-4/5-\delta}). \nonumber
\end{eqnarray}
We now use a Taylor expansion of $E^{\mathcal{X}} V_i(\theta,\mu
_2,x_1)$ with respect to $\theta$. Note that with $A^i = \varepsilon
^i + m_1(X^i_1) - m_1(x_1)$ and $B^i= Y^i -\theta- \mu_2(X_2 ^i) =
\varepsilon^i + m_1(X^i_1) -\theta+ m_2(X^i_2)- \mu_2(X_2 ^i)$
\[
V_i(\theta,\mu_2,x_1) = K_{1,h_1}(x_1,X_1 ^i)\cases{
0, &\quad if $A^i,B^i < 0$,\cr
0, &\quad if $A^i,B^i \geq0$, \cr
B^i, &\quad if $A^i < 0 \leq B^i$, \cr
-B^i, &\quad if $A^i \geq0 > B^i$.}
\]
For $\delta_1,\delta_2>0$ small enough, we get that uniformly for
$|\theta-m_1(x_1)| \leq\delta_1$
\begin{eqnarray*}
&&E^{\mathcal{X}} V_i(\theta,\mu_2,x_1) \\
&&\qquad= \tfrac{1 }{
2}K_{1,h_1}(x_1,X_1 ^i)
f_{\varepsilon|X}(0|X^i) \bigl\{ [m_2(X_2^i)-\mu_2(X_2^i)-\theta
+m_1(x_1) ]^2\\
&&\qquad\quad\hspace*{118.9pt}{} +O_P(n^{-4/5-\delta_2})+O_P\bigl(|\theta-m_1(x_1)|^3\bigr) \bigr\},
\end{eqnarray*}
see (A3). We now apply (\ref{uni1}), (\ref{uni2}) and the fact that
the change of an empirical quantile cannot be larger than the largest
change of an observation. We use these results to analyze the update
$\hat m_1^{\BF}(x_1)$ when we plug\vspace*{2pt} into the iteration formula (\ref
{it:BF}) of the backfitting estimator a choice of $\mu_2=\hat
m_2^{\BF}$ that lies in $\mathcal{M}_n$. By a direct argument, it can
be shown that with probability tending to one the resulting value lies
in an $\delta_1$-neighborhood of $m_1(x_1)$. Thus,\vspace*{1pt} using the above
expansions, we get that, up to terms of order $O_P(n^{-2/5-\delta_3})$
with $\delta_3 > 0$ small enough, the resulting value for the update
$\hat m_1^{\BF}(x_1)$ is equal to the minimum of
\begin{eqnarray*}
&&{\theta\over n} \sum_{i=1}^n K_{1,h_1}(x_1,X_1 ^i) \bigl[I
\bigl(\varepsilon^i + m_1(X^i_1) -
m_1(x_1)\leq0 \bigr) - \alpha\bigr] \\
&&\qquad{} + {1 \over2n} \sum_{i=1}^n K_{1,h_1}(x_1,X_1 ^i)
f_{\varepsilon|X}(0|X^i)
[m_2(X_2^i)-\mu_2(X_2^i)-\theta+m_1(x_1) ]^2.
\end{eqnarray*}
The minimum of this expression is equal to
\begin{eqnarray*}
&&m_1(x_1) - \hat f^w_{X_j}(x_j)^{-1}{1 \over n} \sum
_{i=1}^n K_{1,h_1}(x_1,X_1^i) \bigl[I \bigl(\varepsilon^i + m_1(X^i_1) -
m_1(x_1)\leq0 \bigr) - \alpha\bigr] \\
&&\qquad{} + \hat f^w_{X_j}(x_j)^{-1}{1 \over n} \sum_{i=1}^n
K_{1,h_1}(x_1,X_1 ^i) f_{\varepsilon|X}(0|X^i)
[m_2(X_2^i)-\mu_2(X_2^i) ],
\end{eqnarray*}
where $\hat f^w_{X_j}(x_j)$ has been defined after (\ref
{it:TSBF}). We now use that
\begin{eqnarray*}
&& \hat f^w_{X_j}(x_j)^{-1}{1 \over n} \sum_{i=1}^n
K_{1,h_1}(x_1,X_1 ^i) \bigl[ I \bigl(\varepsilon^i + m_1(X^i_1) -
m_1(x_1)\leq0 \bigr) - I (\varepsilon^i \leq0 ) \bigr] \\
&&\qquad = m_1(x_1) - \hat f^w_{X_j}(x_j)^{-1}{1 \over n}
\sum_{i=1}^n K_{1,h_1}(x_1,X_1 ^i)f_{\varepsilon|X}(0|X^i) m_1(X_1
^i)\\
&&\qquad\quad{} +
O_P(n^{-2/5-\delta})
\end{eqnarray*}
for $\delta> 0$ small enough. This shows that the minimum is equal to
\begin{eqnarray*}
&&\hat f^w_{X_j}(x_j)^{-1}{1 \over n}
\sum_{i=1}^n K_{1,h_1}(x_1,X_1 ^i)f_{\varepsilon|X}(0|X^i) [
m_1(X_1 ^i) + m_2(X_2 ^i) + \eta^i ] \\
&&\qquad{} - \hat f^w_{X_j}(x_j)^{-1}{1 \over n} \sum_{i=1}^n
K_{1,h_1}(x_1,X_1 ^i) f_{\varepsilon|X}(0|X^i)
\mu_2(X_2^i)+O_P(n^{-2/5-\delta}).
\end{eqnarray*}
This expansion holds uniformly for $x_1 \in I_1$ and $\mu_2\in
\mathcal{M}_n$.

To complete the proof, we use the fact that, if one replaces in (\ref
{it:BF}) or (\ref{it:TBF}) the input function $\mu_2 = \hat m_2^{\BF}$
or $\mu_2 = \hat m_2^{*,\BF}$, respectively, by another function that
differs in sup-norm by an amount of order $O_P(n^{-2/5-\Delta_2})$,
then the resulting estimator changes also at most by an amount of order
$O_P(n^{-2/5-\Delta_2})$. In particular, if $\delta< \Delta_2$, this
implies that
\[
{\sup_{a_1+ C_Sn^{-1/5}\leq x_1\leq b_1- C_Sn^{-1/5}}} | \hat
m_1^{\BF}(x_1)-\hat m_1^{*,\BF}(x_1)| =
O_P(n^{-2/5-\delta}).
\]
The other statements of Proposition \ref{TH:2.1} can be proved by using similar
arguments.

%s5.2 ###
\subsection[Proof of Theorem 2.2]{Proof of Theorem \protect\ref{TH:2.2}}
\label{sec:4.2}

We will prove the theorem for the ordinary backfitting estimator. A
proof for the smooth backfitting estimator follows along the same
lines. We only give an outline of the proof. For simplicity, we assume
that the condition (A6) holds with $\rho=1$. Our basic argument runs
as follows. We choose $\hat m_j^{*,\BF,[0]}=\hat m_j^{\BF,[0]}$. By
assumption, these starting values fulfill (A5) and (A6) (with the
choice $\hat m_j^{\BF}=\hat m_j^{*,\BF,[0]}=\hat m_j^{\BF,[0]}$). Thus, we
can apply Proposition \ref{TH:2.1} and we get that the updates $\hat
m_j^{*,\BF,[1]}$ and $\hat m_j^{\BF,[1]}$ fulfill (A7) (with the choices
$\hat m_j^{*,\BF}=\hat m_j^{*,\BF,[1]}$ and $\hat m_j^{\BF}=\hat
m_j^{\BF,[1]}$). We will show below that the updates $\hat
m_j^{*,\BF,[l]}$ of the mean regression\vspace*{1pt} backfitting estimator fulfill
conditions (A5) and (A6) for all $l \geq1$. With this fact, we can use
an iterative argument. Suppose that we know that (A5)--(A7) hold for
$\hat m_j^{*,\BF,[l-1]}$ and $\hat m_j^{\BF,[l-1]}$. Then with our proof
below we get that $\hat m_j^{*,\BF,[l]}$ fulfills (A5) and (A6). By
application of Proposition \ref{TH:2.1}, we get that (A7) holds for $\hat
m_j^{*,\BF,[l]}$ and $\hat m_j^{\BF,[l]}$. Thus, $\hat m_j^{\BF,[l]}$ lies
in a neighborhood of $\hat m_j^{*,\BF,[l]}$ and (A5) and (A6) also hold
for $\hat m_j^{\BF,[l]}$ because they are satisfied by $\hat m_j^{*,\BF,[l]}$.

The bound for the distance between $\hat m_j^{*,\BF,[l]}$ and $\hat
m_j^{\BF,[l]}$ adds up. Each application of Proposition \ref{TH:2.1} adds an
additional term. The additional term increases with $l$. With a careful
analysis of the arguments in the proof of Proposition \ref{TH:2.1}, one gets
that the bounds in (A5) and (A6) have to be multiplied by a factor
$C_*^l$ with a constant $C_* > 1$. If $l \leq C_{\mathrm{iter}} \log n$
with $C_{\mathrm{iter}}>0 $ small enough, we get
%
%e5.5 ###
%
\begin{equation} \label{cl1}
\hat m_j^{\BF,[l]} - \hat m_j^{*,\BF,[l]} = o_P(n^{-2/5}).
\end{equation}
In the second part of the proof, we will show the asymptotic normality
of $\hat m_j^{*,\BF,[C \log n]}$ for $C$ large enough. The minimal
sufficient\vspace*{-1pt} value of $C$ for this result depends on the rate of
convergence of $\hat m_j^{*,\BF,[0]}$ to $m_j$. If this rate is
$n^{-2/5}$, then it can be made as small as one likes. For slower
rates, one needs larger values of $C$. If the rate is fast enough, one
can choose $C< C_{\mathrm{iter}}$. In this case, we can apply (\ref{cl1})
and we get the same asymptotic normality result for $\hat
m_j^{*,\BF,[C_{\mathrm{iter}} \log n]}$. This will conclude the proof of
Theorem \ref{TH:2.2}.

We now prove that the updates $\hat m_j^{*,\BF,[l]}$ fulfill the
conditions (A5) and (A6) for all $l \geq1$. For this purpose, we rewrite
(\ref{it:TBF}) as
%
%e5.6 ###
%
\begin{eqnarray} \label{star2}
&&\hat m_j^{*,\BF,[l]} (x_j) - m_j(x_j) \nonumber\\
&&\qquad = \tilde m_j^{*,A} (x_j)+ \tilde m_j^{*,B} (x_j) + \tilde
m_j^{*,C, [l]} (x_j)- \hat m_0^{*,\BF} \\
&&\qquad\quad{} - \sum_{k=1, \not= j}^d \int \bigl[ \hat
m_k^{*,\BF,[l_{k,j}]} (x_k) - m_k (x_k) \bigr]
f_{X_k|X_j}^{n,w}(x_k|x_j) \,dx_k, \nonumber
\end{eqnarray}
where $l_{k,j}= l+1$ for $k< j$, $l_{k,j}=l$ for $k > j$, and
\begin{eqnarray*}
\tilde m_j^{*,A} (x_j)&=& { n^{-1} \sum_{i=1}^n f_{\varepsilon
|X}(0|X^i) K_{j,h_j}(x_j,X_j^i) \eta^i \over\hat f_{X_j}^w(x_j) },\\
\tilde m_j^{*,B} (x_j) &=& { n^{-1}\sum_{i=1}^n f_{\varepsilon
|X}(0|X^i) K_{j,h_j}(x_j,X_j^i) [m_j(X_j^i)-m_j(x_j)] \over\hat
f_{X_j}^w(x_j) }, \\
\tilde m_j^{*,C, [l]} (x_j)
&=& - \sum_{k=1, \not= j}^d\Biggl( n^{-1} \sum_{i=1}^n f_{\varepsilon
|X}(0|X^i) K_{j,h_j}(x_j,X_j^i) \\
&&\qquad\hspace*{56.2pt}{}\times\bigl[\hat m_k^{*,\BF,[l_{k,j}]}
(X_k^i)-m_k(X_k^i)\bigr] \Biggr)(\hat f_{X_j}^w(x_j) )^{-1}\\
&&{} + \sum_{k=1, \not= j}^d \int{ \bigl[ \hat
m_k^{*,\BF,[l_{k,j}]} (x_k) - m_k (x_k) \bigr]
f_{X_k|X_j}^{n,w}(x_k|x_j) \,dx_k},\\
f_{X_k|X_j}^{n,w}(u_k|x_j)&=&{\int f_{\varepsilon|X}(0|u)
K_{j,h_j}(x_j,u_j) f_X(u) \,du_{-k} \over\int f_{\varepsilon|X}(0|v)
K_{j,h_j}(x_j,v_j) f_X(v) \,dv}.
\end{eqnarray*}
The iteration\vspace*{1pt} (\ref{star2}) can be analyzed as the smooth backfitting
algorithm in \citet{MLN99}. With $\hat
m_+^{*,\BF,[l]}(x) = \hat m_1^{*,\BF,[l]}(x_1) + \cdots+ \hat
m_d^{*,\BF,[l]}(x_d)$ and $m_+(x) = m_1(x_1) + \cdots+ m_d(x_d) $, we
can write a full cycle of iterations (\ref{star2}) as
%
%e5.7 ###
%
\begin{eqnarray}\label{plus}
\hat m_+^{*,\BF,[l+1]} - m_+
&=& \tilde m_{\oplus}^{*,A} + \tilde m_{\oplus}^{*,B} +
\tilde m_{\oplus}^{*,C, [l]} - \hat m_0^{*,\BF}\nonumber\\[-8pt]\\[-8pt]
&&{}+ T_{n,+} \bigl( \hat
m_+^{*,\BF,[l]} - m_+ - \mu_l \bigr)+ \mu_l, \nonumber
\end{eqnarray}
where $\tilde m_{\oplus}^{*,A}$, $ \tilde m_{\oplus}^{*,B}$ and $
\tilde m_{\oplus}^{*,C, [l]}$ are some functions, $T_{n,+}$ is an
operator that acts on additive mean zero functions in
$L_2(f_{\varepsilon|X}(0|\cdot)f_X(\cdot))$, and $\mu_l= \int(\hat
m_+^{*,\BF,[l]} - m_+)(x) f_X(x) f_{\varepsilon|X}(0|x)\,dx$. We used
${\oplus}$ (not $+$) as subindex in $\tilde m_{\oplus}^{*,A}$ because
it is not the sum of $\tilde m_{j}^{*,A}$. The operator $T_{n,+}$
converges to an operator $T_{+}$ that is based on an iterative
application of the linear transformations for the additive components
$g_j$ of an additive function $g_+$:
\[
g_j \to - \sum_{k=1, \not= j}^d \int{ g_k (x_k)
f_{X_k|X_j}^w(x_k|x_j) \,dx_k}.
\]
More precisely, the kernel function of $T_{n,+}$ converges to the
kernel function of $T_{+}$, with respect to the sup-norm.

Arguing as in the proof of Lemma 1 in \citet{MLN99},
one can show that $T_+$ is a positive self-adjoint operator with
operator norm strictly less than one, $\|T_{+}\| < 1$, and with $\|T_j
m \|_{\infty} \leq D \|m\|_2$ for a constant $D> 0$. Here, $T_jm$ is
the $j$th additive component of $T_+m$. This gives with constants $0 <
D^{\prime} < 1$ and $D^{\prime\prime}>0$ for $n$ large enough
%
%e5.8 ###
%
\begin{equation} \label{star3}
\|T_{n,+}\| < D^{\prime}.
\end{equation}
Furthermore, we have
%
%e5.9 ###
%
\begin{equation} \label{star4}
\|T_{n,j} m \|_{\infty} \leq D^{\prime\prime} \|m\|_2,
\end{equation}
where $T_{n,j}m$ is the $j$th additive component of $T_{n,+}m$.
Iterative application of (\ref{plus}) gives
\[
\hat m_+^{*,\BF,[l]} - m_+ = \hat m_+^{*,A,[l]} + \hat m_+^{*,B,[l]} +
\hat m_+^{*,C, [l]} - \hat m_0^{*,\BF}+ \bar T_{n,+}^l \bigl( \hat
m_+^{*,\BF,[0]} - m_+ \bigr),
\]
where $\bar T_{n,+}$ is an extension of $T_{n,+}$ to a nonzero mean
function by putting
$\bar T_{n,+} g = T_{n,+}(g- \mu_g) + \mu_g$ with $\mu_g=\int g(x)
f_X(x) f_{\varepsilon|X}(0|x)\,dx$, and
\begin{eqnarray*}
\hat m_+^{*,A,[l]} &=& \sum_{r=0}^{l-1} \bar T_{n,+}^r \tilde
m_{\oplus}^{*,A}, \\
\hat m_+^{*,B,[l]} &=& \sum_{r=0}^{l-1} \bar T_{n,+}^r \tilde
m_{\oplus}^{*,B}, \\
\hat m_+^{*,C,[l]} &=& \sum_{r=0}^{l-1} \bar T_{n,+}^{l-r-1} \tilde
m_{\oplus}^{*,C, [r]}.
\end{eqnarray*}
Using standard bounds on $\tilde m_j^{*,A}$ and $\tilde m_j^{*,B}$, it
can be verified that
%
%e5.12 ###
%e5.11 ###
%e5.10 ###
%
\begin{eqnarray}
\label{eq:a1}
\sup_{x_j \in I_j, l \geq1} \bigl|\hat m_j^{*,A,[l]}(x_j)\bigr| &=&
O_P(n^{-2/5}), \\
\label{eq:a2}
\sup_{x_j \in I_j, l \geq1} \bigl|\hat m_j^{*,B,[l]}(x_j)\bigr| &=&
O_P(n^{-1/5}), \\
\label{eq:a3}
\sup_{a_j + C_S h_j \leq x_j \leq b_j - C_S h_j, l \geq1} \bigl|\hat
m_j^{*,B,[l]}(x_j)\bigr| &=& O_P(n^{-2/5}),
\end{eqnarray}
where for an additive function $g_{+} $ we denote by $g_j$ its $j$th
additive component.

We now argue that for a constant $C_T>0$
%
%e5.13 ###
%
\begin{equation}
\label{eq:a4}
\sup_{x_j \in I_j, l \geq1} \bigl|\bar T_{n,j} \bar T_{n,+}^{l-1}
\bigl( \hat m_j^{*,\BF,[0]} - m_j \bigr)(x_j) \bigr| \leq C_T
\kappa_n,
\end{equation}
where
\begin{eqnarray*}
\kappa_n &=& \sup_{1 \leq j \leq d} \Bigl[\sup_{a_j + C_S h_j \leq
x_j \leq b_j - C_S h_j} \bigl| \hat m_j^{*,\BF,[0]} - m_j \bigr| (x_j)\\
&&\hspace*{28.8pt}\hspace*{13.4pt}{} +
n^{-1/5} \sup_{a_j \leq x_j \leq b_j } \bigl| \hat m_j^{*,\BF,[0]} -
m_j \bigr| (x_j) \Bigr].
\end{eqnarray*}
For a proof of this claim, one applies (\ref{star3}) and (\ref
{star4}). Also, we argue
that
%
%e5.14 ###
%
\begin{equation}
\label{eq:a5}
\sup_{x_j \in I_, l \geq1} \bigl|\hat m_j^{C,[l]}(x_j)\bigr| = o_P(n^{-2/5}).
\end{equation}
For a proof of (\ref{eq:a5}), we note that
\[
\sup_{x_j \in I_j, l \geq1} \bigl|\tilde m_j^{C,[l]}(x_j)\bigr| = o_P(n^{-2/5}).
\]
This follows by empirical process theory. One uses the fact that $\hat
m_k^{*,\BF,[l-1]}-m_k$ lies in a class of functions that have second
derivatives absolutely bounded by $C_{\xi} n ^{\xi}$ with $\xi>0$
being arbitrarily small and constant $C_\xi$ depending on $\xi$. This
can be shown by using that the same bound applies for $\tilde m_j^A$
and $\tilde m_j^B$, and that the kernels of the operators $T_+$ and
$T_j$ have an absolutely bounded second derivative [see (A9)], and then
applying an iterative argument.

The bounds\vspace*{1pt} at (\ref{eq:a1})--(\ref{eq:a5}) imply that $\hat
m_j^{*,\BF,[l]}$ fulfills (A5) uniformly for \mbox{$l \geq1$}. Using the
smoothness considerations in the previous paragraph, we get that $\hat
m_j^{*,\BF,[l]}$ fulfills\vspace*{1pt} (A6) uniformly for $l \geq1$. Thus, we get by
an iterative application of Proposition \ref{TH:2.1} that (\ref{cl1}) holds.

It remains to show the asymptotic normality result for $\hat m_j^{\BF,
\mathrm{iter}}=\break\hat m_j^{\BF,[C_{\mathrm{iter}} \log n]}$ with
$C_{\mathrm{iter}}$ large enough. Using the above arguments, we have for
$C_{\mathrm{iter}}$ large enough that
\[
\hat m_j^{*,\BF,\mathrm{iter}}(x_j) - m_j(x_j) = \hat m_j^{A,
[C_{\mathrm{iter}} \log n]} (x_j)+ \hat m_j^{B,[C_{\mathrm{iter}} \log n]} (x_j)+
o_P(n ^{-2/5}).
\]
We argue that
%
%e5.16 ###
%e5.15 ###
%
\begin{eqnarray}
\label{asnor1}
\sup_{l\geq1} \bigl| \hat m_j^{A,[l]} (x_j)- \tilde m_j^{A} (x_j)\bigr| &=&
o_P(n^{-2/5}), \\
\label{asnor2}
h_j^{-2} \hat m_j^{B,[l]} (x_j) &\to& \beta_j(x_j) \qquad\mbox{as }
l\to\infty.
\end{eqnarray}
These two claims imply that
\[
\hat m_j^{*,\BF}(x_j) - m_j(x_j) = \tilde m_j^{A} (x_j)+ h_j^2 \beta
_j(x_j)+ o_P(n ^{-2/5}).
\]
This expansion shows the desired asymptotic limit result by using a
standard smoothing limit result for $\tilde m_j^{A} (x_j)$.\vspace*{1pt}

We prove (\ref{asnor1}) and (\ref{asnor2}). Claim (\ref{asnor1})
follows from standard smoothing theory as in
\citet{MLN99}. For a proof of (\ref{asnor2}), we define
$\beta_j^{[l]}(x_j)= \beta_j^*(x_j) - \sum_{k=1,\not=j}^d \int
\beta_k^{[l_{k,j}]}(x_k) f^w_{X_k|X_j}(x_k|x_j) \,dx_k$ with $\beta
_j^{[0]}(x_j) \equiv0$. Similarly, as in (\ref{plus}), we can write a
full cycle of these iterations as
%
%e5.17 ###
%
\begin{equation} \label{foc}
\beta_+^{[l+1]} = \beta_{\oplus}^* + \bar T_+ \beta_+^{[l]},
\end{equation}
where $\beta_{\oplus}^*$ is some additive function, $\beta
_+^{[l]}(x)$ is equal to $\beta_1^{[l]}(x_1)+ \cdots+\beta
_d^{[l]}(x_d)$ and $\bar T_+$ is an extension of $ T_+$ defined by
$\bar T_+g= T_+(g- \mu_g) + \mu_g$ with $\mu_g$ defined as above.
Note that we get $\beta_+^{[l]}= \sum_{r=0}^{l-1} \bar T_{+}^r \beta
_{\oplus}^*$. This expansion shows that
%
%e5.18 ###
%
\begin{eqnarray}\label{adda1}
\sup_{x_j \in I_j, l \geq1} \bigl|\hat
m_j^{*,B,[l]}(x_j)-h_j^{2}\beta_j^{[l]} \bigr| &=& o_P(n^{-1/5}),
\nonumber\\[-8pt]\\[-8pt]
\sup_{a_j + C_S h_j \leq x_j \leq b_j - C_S h_j, l
\geq1} \bigl|\hat m_j^{*,B,[l]}(x_j)-h_j^{2}\beta_j^{[l]} \bigr| &=&
o_P(n^{-2/5}).\nonumber
\end{eqnarray}
Furthermore, we get that the term $\beta_+^{[l]} - \sum_{j=1}^d
[h_j^{-2} m_j^{\prime} (x_j) \int(u_j-x_j) K_{j,h_j}(x_j,\break u_j) \,du_j
-\mu_{2,K} {1 \over2} m_j^{\prime\prime}(x_j) ]$ converges to
$\mu_{2,K} \beta_+^{**}$ as $l \to\infty$, where $(\beta
_1^{**},\ldots,\beta_d^{**})$ is the minimizer of
\[
\int\Biggl[ \sum_{j=1}^d \biggl(m_j^{\prime}
(x_j) {{\partial/\partial
x_j}f_{\varepsilon, X}(0,x) \over f_{\varepsilon, X}(0,x) } - \beta
_j^{**}(x_j) \biggr) \Biggr] ^2 f_{\varepsilon, X}(0,x) \,dx.
\]
This follows because the updating (\ref{foc}) is given by the
first-order conditions of this minimization problem. Together with
(\ref{adda1}), this implies (\ref{asnor2}).

\printaddresses


\begin{thebibliography}{99}

%b1 ###
\bibitem[\protect\citeauthoryear{Bahadur}{1966}]{B66}
\textsc{Bahadur, R. R.} (1966). A
note on quantiles in large samples.
\textit{Ann. Math. Statist.} \textbf{37} 577--580.
\MR{0189095}

%b2 ###
\bibitem[\protect\citeauthoryear{Buja, Hastie and Tibshirani}{1989}]{BHT89}
\textsc{Buja, A., Hastie, T.} and \textsc{Tibshirani, R.} (1989). Linear
smoothers and additive models (with discussion).
\textit{Ann. Statist.} \textbf{17} 453--555.
\MR{0994249}

%b3 ###
\bibitem[\protect\citeauthoryear{Chaudhuri}{1991}]{C91}
\textsc{Chaudhuri, P.} (1991).
Nonparametric estimates of regression quantiles and their
Bahadur representation. \textit{Ann. Statist.} \textbf{19} 760--777.
\MR{1105843}

%b4 ###
\bibitem[\protect\citeauthoryear{Fan and Gijbels}{1996}]{FG96}
\textsc{Fan, J.} and \textsc{Gijbels, I.}
(1996). \textit{Local Polynomial Modelling and Its Applications}.
Chapman and Hall, London.
\MR{1383587}

%b5 ###
\bibitem[\protect\citeauthoryear{Furno}{2004}]{F04}
\textsc{Furno, M.} (2004). ARCH tests
and quantile regressions. \textit{J. Stat. Comput. Simul.}
\textbf{74} 277--292.
\MR{2059314}

%b6 ###
\bibitem[\protect\citeauthoryear{Horowitz and Lee}{2005}]{HL05}
\textsc{Horowitz, J.} and \textsc{Lee, S.}
(2005). Nonparametric estimation of an additive
quantile regression model. \textit{J. Amer. Statist. Assoc.} \textbf{100}
1238--1249.
\MR{2236438}

%b7 ###
\bibitem[\protect\citeauthoryear{Horowitz and Mammen}{2004}]{HM04}
\textsc{Horowitz, J.} and \textsc{Mammen, E.}
(2004). Nonparametric estimation of an additive
model with a link function. \textit{Ann. Statist.} \textbf{32} 2412--2443.
\MR{2153990}

%b8 ###
\bibitem[\protect\citeauthoryear{Jones and Hall}{1990}]{JH90}
\textsc{Jones, M. C.} and \textsc{Hall, P.}
(1990). Mean squared error properties of kernel estimates
of regression quantiles. \textit{Statist. Probab. Lett.} \textbf{10}
283--289.
\MR{1069903}

%b9 ###
\bibitem[\protect\citeauthoryear{Lee, Lee and Park}{2006}]{LLP06}
\textsc{Lee, Y. K., Lee, E. R.} and \textsc{Park, B. U.} (2006).
Conditional quantile estimation by
local logistic regression. \textit{J. Nonparametr. Stat.}
\textbf{18} 357--373.
\MR{2284188}

%b10 ###
\bibitem[\protect\citeauthoryear{Linton and Nielsen}{1995}]{LN95}
\textsc{Linton, O.} and \textsc{Nielsen, J. P.} (1995). A kernel method of
estimating structured nonparametric regression based on marginal
integration. \textit{Biometrika} \textbf{82} 93--101.
\MR{1332841}

%b11 ###
\bibitem[\protect\citeauthoryear{Lu and Yu}{2004}]{LY04}
\textsc{Lu, Z.} and \textsc{Yu, K.} (2004).
Local linear additive quantile regression.
\textit{Scand. J. Statist.} \textbf{31} 333--346.
\MR{2087829}

%b12 ###
\bibitem[\protect\citeauthoryear{Mammen, Linton and Nielsen}{1999}]{MLN99}
\textsc{Mammen, E., Linton, O.} and \textsc{Nielsen, J. P.} (1999). The
existence and asymptotic properties of a backfitting projection
algorithm under weak conditions. \textit{Ann. Statist.} \textbf{27}
1443--1490.
\MR{1742496}

%b13 ###
\bibitem[\protect\citeauthoryear{Mammen and Park}{2006}]{MP06}
\textsc{Mammen, E.} and \textsc{Park, B. U.}
(2006). A simple smooth backfitting method for additive models.
\textit{Ann. Statist.} \textbf{34} 2252--2271.
\MR{2291499}

%b14 ###
\bibitem[\protect\citeauthoryear{Opsomer}{2000}]{O00}
\textsc{Opsomer, J. D.} (2000).
Asymptotic properties of backfitting
estimators. \textit{J. Multivariate Anal.} \textbf{73}
166--179.
\MR{1763322}

%b15 ###
\bibitem[\protect\citeauthoryear{Opsomer and Ruppert}{1997}]{OR97}
\textsc{Opsomer, J. D.} and \textsc{Ruppert, D.} (1997). Fitting a bivariate
additive model by local polynomial regression.
\textit{Ann. Statist.} \textbf{25} 186--211.
\MR{1429922}

%b16 ###
\bibitem[\protect\citeauthoryear{Yu and Jones}{1998}]{YJ98}
\textsc{Yu, K.} and \textsc{Jones, M. C.}
(1998). Local linear quantile regression.
\textit{J. Amer. Statist. Assoc.} \textbf{93} 228--237.
\MR{1614628}

%b17 ###
\bibitem[\protect\citeauthoryear{Yu, Lu and Stander}{2003}]{YLS03}
\textsc{Yu, K., Lu, Z.} and \textsc{Stander, J.} (2003).
Quantile regression: Applications and
current research areas. \textit{The Statistician} \textbf{52} 331--350.
\MR{2011179}

%b18 ###
\bibitem[\protect\citeauthoryear{Yu, Park and Mammen}{2008}]{YPM08}
\textsc{Yu, K., Park, B. U.} and \textsc{Mammen, E.} (2008).
Smooth backfitting in generalized
additive models. \textit{Ann. Statist.} \textbf{36} 228--260.
\MR{2387970}

\end{thebibliography}
\end{document}

\begin{frontmatter}

\title{Correction\\
Backfitting and smooth backfitting for additive quantile models}\vspace*{4pt}
\runtitle{Correction}
\pdftitle{Correction
Backfitting and smooth backfitting for additive quantile models}

\textit{Ann. Statist.} \textbf{38} (2010) 2857--2883

\begin{aug}
\author[A]{\fnms{Young K.} \snm{Lee}\ead[label=e1]{youngklee@kangwon.ac.kr}},
\author[B]{\fnms{Enno} \snm{Mammen}\ead[label=e2]{emammen@rumms.uni-mannheim.de}}
\and
\author[C]{\fnms{Byeong U.} \snm{Park}\corref{}\ead[label=e3]{bupark@stats.snu.ac.kr}}
\runauthor{Y. K. Lee, E. Mammen and B. U. Park}
\affiliation{Kangwon National University,
University of Mannheim and\break Seoul National University}
\address[A]{Y. K. Lee\\
Department of Statistics\\
Kangwon National University\\
Chuncheon 200-701\\
Korea\\
\printead{e1}}
\address[B]{E. Mammen\\
Department of Economics\\
University of Mannheim\\
68131 Mannheim, L7, 3-5\\
Germany\\
\printead{e2}}
\address[C]{B. U. Park\\
Department of Statistics\\
Seoul National University\\
Seoul 151-747\\
Korea\\
\printead{e3}} %adresu isvedimo komanda gale!
\end{aug}

% HISTORY:
\received{\smonth{9} \syear{2012}}

% ABSTRACT

% KEYWORDS
% Pirmas kwd is didziosios raides

\end{frontmatter}

In Theorem 2.2 on page 2865 of \cite{lee2010} we wrongly stated that
the asymptotic biases
of the ordinary and of the smooth backfitting estimator were the same.
In fact, the bias formulas for the two methods are different. The
theorem should be modified as follows.

\renewcommand{\thetheorem}{2.2}
\begin{theorem}\label{theorem2.2}
Let $\alpha_j(u)=m_j'(u)\int(v-u)K_{j,h_j}(u,v) \,dv /\int
K_{j,h_j}(u,v) \,dv$ and $\mu_{2,K}=\int v^2 K(v) \,dv$. Assume that
\textup{(A1)--(A4), (A8)} and \textup{(A9)} hold, that \textup{(A5)}
and \textup{(A6)} are satisfied by $\hat m_j^{\mathrm{BF}}=\hat
m_j^{\mathrm{BF},[0]}$ and $\hat m_j^{\mathrm{SBF}}=\hat
m_j^{\mathrm{SBF},[0]}$ ($j=2,\ldots,d$) with $\xi,\Delta_2,\Delta_3,
{2 \over5} - {1+ \rho\over2 + 3 \rho}{4 \over5} - \Delta_1 >0$ small
enough,\vspace*{2pt} and that $w_j(a_j+x(b_j-a_j)) \le C x(1-x)$ for
all $0 \le x \le1$ and for some positive constant $C$.
Then,\vspace*{1pt} we get for $\hat m_j^{l,{\mathrm{iter}}}=\hat
m_j^{l,[C_{\mathrm{iter}} \log n]}$ with\vspace*{1pt} an appropriate
choice of $C_{\mathrm{iter}}=C_{\mathrm{iter},l}$ ($l= \mathrm{BF}$ and
$l=\mathrm{SBF}$) that for $a_j < x_j < b_j$
\begin{eqnarray*}
&&\sqrt{nh_j} \bigl[\hat m_j^{l,\mathrm{iter}}
(x_j) -m_j(x_j) - \beta_j^l(x_j)
\bigr]
\\
&&\qquad \to N \biggl(0,\frac{\alpha(1-\alpha)}{f_{\varepsilon
,X_j}(0,x_j)^2}f_{X_j}(x_j) \int
K^2(u) \,du \biggr)
\end{eqnarray*}
in distribution. Here, $\beta_j^l(x_j) = \beta_j^{*,l}(x_j) - \int
\beta_j^{*,l}(u_j)w_j(u_j) \,du_j$, and
$(\beta_j^{*,l}\dvtx  1 \le j \le d)$ for $l=\mathrm{BF}$ is the solution of
the system of integral equations
\begin{eqnarray}
0&=&\int\Biggl[\alpha_j(x_j) +h_j^2
\mu_{2,K} m_j'(x_j)\,
\frac{\partial f_{\varepsilon, X}(0,x)/\partial x_j} {
f_{\varepsilon, X}(0,x)}\nonumber\\
&&\hspace*{49pt}{} + \frac{1}{2}h_j^2
\mu_{2,K} m_j''(x_j)
 - \sum_{k=1}^d \beta_k^{*,l}(x_k)
\Biggr]f_{\varepsilon,X}(0,x) \,dx_{-j},\nonumber\\
\eqntext{1 \le j \le d,}
\end{eqnarray}
while, for $l=\mathrm{SBF}$, $\beta_j^{*,l}(x_j)=\frac{1}{2}h_j^2 \mu
_{2,K} m_j''(x_j) + \mu_{2,K}\beta_j^{**}(x_j)$ where $(\beta_j^{**}\dvtx  1
\le j \le d)$ is a tuple of functions that minimizes
\[
\int\Biggl[\sum_{j=1}^d
\biggl(h_j^2 m_j'(x_j)\,
\frac{\partial
f_{\varepsilon,X}(0,x)/\partial x_j} {
f_{\varepsilon,X}(0,x)}- \beta_j^{**}(x_j) \biggr)
\Biggr]^2 f_{\varepsilon, X}(0,x) \,dx.
\]
\end{theorem}

In addition, the assumption (A9) on page 2864 should be modified as follows:
\begin{longlist}[(A9)]
\item[(A9)] The functions $f_{X_j,X_k}^w(x_j,x_k)$ and
$f_{X_j}^w(x_j)$ have second derivatives w.r.t. $x_j$ that are bounded
over $x_j\in I_j$, $x_k\in I_k$ for all $1\leq j \neq k \leq d$. Also,
$f_{X_j,X_k}^w(x_j,x_k)$ and its first derivative w.r.t. $x_j$ is
continuous in $x_k$ for all $1\leq j \neq k \leq d$.
\end{longlist}

The corrected statement of Theorem 2.2 and the modified assumption (A9)
do not require a different proof than what is given in \cite{lee2010}.

%suskaldyti doi

% imsref loaded by lrinkeviciute, 2012-11-28 11:17:30

\printaddresses